\documentclass[10pt]{elsart}
\usepackage{graphicx}
\usepackage{amsmath}
\usepackage{amsfonts}
\usepackage{enumerate}
\usepackage{latexsym}
\usepackage{epsfig}

\newtheorem{theorem}{Theorem}[section]
\newtheorem{lemma}[theorem]{Lemma}
\newtheorem{coroll}[theorem]{Corollary}

\def\proofbox{\begin{picture}(6.5,6.5)
\put(0,0){\framebox(6.5,6.5){}}\end{picture}}
\newenvironment{proof}{\noindent{\it Proof.\quad}}{\hfill\proofbox}

\def\rank{\mathrm{rank}\;}

\begin{document}
\begin{frontmatter}

\title{Superinjective Simplicial Maps of Complexes of Curves and Injective
Homomorphisms of Subgroups of Mapping Class Groups }
\author{Elmas Irmak}

\maketitle

%pagestyle{fancy}
\renewcommand{\sectionmark}[1]{\markright{\thesection. #1}}

%\fancyhf{}
%\fancyhead[LE]{{\sc \thepage\hfill{On Injective Homomorphisms
%of Mapping Class Groups}}}
%\fancyhead[RO]{{\sc \rightmark\hfill\thepage}}
%\fancyfoot[LE]{{\tiny Michigan State University -- Department of
%    Mathematics}}
%\fancyfoot[RO]{{\tiny Elmas Irmak}}

\thispagestyle{empty}
\maketitle

\address{Department of Mathematics, University of Michigan,\\
Ann Arbor, MI 48109, USA}

%      e-mail: irmak@math.msu.edu }\\ [1.5cm]

\begin{abstract}
Let $S$ be a closed, connected, orientable surface of genus at
least 3, $\mathcal{C}(S)$ be the complex of curves on $S$ and
$Mod_S^*$ be the extended mapping class group of $S$. We prove
that a simplicial map, $\lambda: \mathcal{C}(S) \rightarrow
\mathcal{C}(S)$, preserves nondisjointness (i.e. if $\alpha$ and
$\beta$ are two vertices in $\mathcal{C}(S)$ and $i(\alpha, \beta)
\neq 0$, then $i(\lambda(\alpha), \lambda(\beta)) \neq 0$) iff it
is induced by a homeomorphism of $S$. As a corollary, we prove
that if $K$ is a finite index subgroup of $Mod_S^*$ and $f : K
\rightarrow Mod_S^*$ is an injective homomorphism, then $f$ is
induced by a homeomorphism of $S$ and $f$ has a unique extension
to an automorphism of $Mod_S^*$.
\end{abstract}
\vspace{0.12in}


\end{frontmatter}

\section{Introduction} Let $R$ be a compact orientable surface
possibly having nonempty boundary. The mapping class group,
$Mod_R$, of $R$ is the group of isotopy classes of orientation
preserving homeomorphisms of $R$. The extended mapping class
group, $Mod_R^*$, of $R$ is the group of isotopy classes of all
(including orientation reversing) homeomorphisms of $R$.\\
\indent In section 2, we give the notations and the terminology
that we use in the paper.\\
\indent In section 3, we introduce the notion of superinjective
simplicial maps of the complex of curves, $\mathcal{C}(S)$, of a
closed, connected, orientable surface $S$ and prove some
properties of these maps.\\
\indent In section 4, we prove that a superinjective simplicial
map, $\lambda$, of the complex of curves induces an injective
simplicial map on the complex of arcs, $\mathcal{B}(S_c)$, where
$c$ is a nonseparating circle on $S$ and we prove that $\lambda$
is induced by a homeomorphism of $S$.\\
\indent In section 5, we prove that if $K$ is a finite index
subgroup of $Mod_S^*$ and $f:K \rightarrow Mod_S^*$ is an
injective homomorphism, then f induces a superinjective simplicial
map of $\mathcal{C}(S)$ and $f$ is induced by a homeomorphism of
$S$.

Our main results are motivated by the following theorems of Ivanov
and the theorem of Ivanov and McCarthy.

\begin{theorem} {\bf (Ivanov) \cite{Iv1}}
\label{theoremA} Let $S$ be a compact, orientable surface possibly with nonempty boundary.
Suppose that the genus of $S$ is at least 2. Then, all automorphisms of $ \mathcal{C}(S)$
are given by elements of $Mod_S^*$. More precisely, if $S$ is not a closed surface of genus 2,
then $Aut(\mathcal{C}(S))= Mod_S ^*$. If $S$ is a closed surface of genus 2,
then $Aut(\mathcal{C}(S))= Mod_S ^* /Center(Mod_S ^*)$.
\end{theorem}

\begin{theorem} {\bf (Ivanov) \cite{Iv1}}
\label{theoremB} Let $S$ be a compact, orientable surface possibly with nonempty boundary.
Suppose that the genus of $S$ is at least 2 and $S$ is not a closed surface of genus 2.
Let $\Gamma_1, \Gamma_2$ be finite index subgroups of $Mod_S ^*$. Then, all isomorphisms
$\Gamma_1 \rightarrow \Gamma_2$ have the form $x \rightarrow gxg^{-1}, g \in Mod_S^*$.
\end{theorem}

Theorem 1.1 and Theorem 1.2 were extended to all surfaces of genus
0 and 1 with the exception of spheres with $\leq$ 4 holes and tori
with $\leq$ 2 holes by M.Korkmaz. These extensions were also
obtained by F.Luo independently.

\begin{theorem} {\bf (Ivanov, McCarthy) \cite{IMc}}
\label{theoremC}
Let $S$ and $S'$ be compact, connected, orientable surfaces.  Suppose that the genus of $S$ is at least 2,
$S'$ is not a closed surface of genus 2, and the maxima of ranks of abelian subgroups of $Mod_S$ and
$Mod_{S'}$ differ by at most one. If $h: Mod_S \rightarrow Mod_{S'}$ is an injective
homomorphism, then  $h$ is induced by a homeomorphism $H: S \rightarrow S'$,
(i.e. $h([G])=[HGH^{-1}]$, for every orientation preserving homeomorphism $G : S \rightarrow S$).
\end{theorem}

The main results of the paper are the following two theorems:

\begin{theorem}
\label{theorem1} Let $S$ be a closed, connected, orientable surface of
genus at least 3. A simplicial map, $\lambda : \mathcal{C}(S) \rightarrow \mathcal{C}(S)$,
is superinjective if and only if $\lambda$ is induced by a homeomorphism of $S$.
\end{theorem}

\begin{theorem}
\label{theorem2} Let $S$ be a closed, connected, orientable
surface of genus at least 3. Let $K$ be a finite index subgroup of
$Mod_S^*$ and $f$ be an injective homomorphism, $f:K \rightarrow
Mod_S^*$. Then $f$ is induced by a homeomorphism of the surface
$S$ (i.e. for some $g \in Mod_S^*$, $f(k) = gkg^{-1}$ for all $k
\in K$) and $f$ has a unique extension to an automorphism of
$Mod_S^*$.
\end{theorem}

Theorem \ref{theorem2} is deduced from Theorem \ref{theorem1}.
Theorem \ref{theorem1} generalizes the closed case of Ivanov's
Theorem \ref{theoremA} for surfaces of genus at least 3. Theorem
\ref{theorem2} generalizes Ivanov's Theorem \ref{theoremB} in the
case of closed surfaces and Ivanov and McCarthy's Theorem
\ref{theoremC} in the case when the surfaces are the same and
closed. In our proof, some advanced homotopy results about
$\mathcal{C}(S)$ used by Ivanov are replaced by elementary surface
topology arguments.

\section{Notations and Preliminaries}
A \textit{circle} on a surface, $R$, of genus $g$ with $b$
boundary components is a properly embedded image of an embedding
$S^{1} \rightarrow R$. A circle on $R$ is said to be
\textit{nontrivial} (or \textit{essential}) if it doesn't bound a
disk and it is not homotopic to a boundary component of $R$. Let C
be a collection of pairwise disjoint circles on $R$. The surface
obtained from $R$ by cutting along $C$ is denoted by $R_C$. Two
circles $a$ and $b$ on $R$ are called \textit{topologically
equivalent} if there exists a homeomorphism $F: R \rightarrow R$
such that $F(a)=b$. The isotopy class of a Dehn twist about a
circle $a$, is denoted by $t_{\alpha}$, where $[a] = \alpha$.

Let $\mathcal{A}$ denote the set of isotopy classes of nontrivial
circles on $R$. The \textit{geometric intersection number}
$i(\alpha, \beta)$ of $\alpha, \beta \in \mathcal{A}$ is the
minimum number of points of $a \cap b$ where $a \in \alpha$ and $b
\in \beta$.

A mapping class, $g \in Mod_R^*$, is called \textit{pseudo-Anosov}
if $\mathcal{A}$ is nonempty and if $g ^n (\alpha) \neq \alpha$,
for all $\alpha$ in $\mathcal{A}$ and any $n \neq 0$. $g$ is
called \textit{reducible} if there is a nonempty subset $
\mathcal{B} \subseteq \mathcal{A}$ such that a set of disjoint
representatives can be chosen for $\mathcal{B}$ and
$g(\mathcal{B}) = \mathcal{B}$. In this case, $ \mathcal{B}$ is
called a \textit{reduction system} for $g$. Each element of
$\mathcal{B}$ is called a \textit{reduction class} for $g$. A
reduction class, $\alpha$, for $g$, is called an \textit{essential
reduction class} for $g$, if for each $\beta \in \mathcal{A}$ such
that $i(\alpha, \beta) \neq 0$ and for each integer $m \neq 0$,
$g^m (\beta) \neq \beta$. The set, $\mathcal{B}_g$, of all
essential reduction classes for $g$ is called the
\textit{canonical reduction system} for $g$. The correspondence $g
\rightarrow \mathcal{B}_g$ is canonical. In particular, it
satisfies $g(\mathcal{B}_h) = \mathcal{B}_{ghg^{-1}}$ for all $g,
h$ in $Mod_R^*$.

The \textit{complex of curves}, $\mathcal{C}(R)$, on $R$ is an
abstract simplicial complex, as given in \cite{Sp}, with vertex
set $\mathcal{A}$ such that a set of $n$ vertices $\{{
\alpha_{1}}, {\alpha_{2}}, ..., {\alpha_{n}}\}$ forms an $n-1$
simplex if and only if ${\alpha_{1}}, {\alpha_{2}},...,
{\alpha_{n}}$ have pairwise disjoint representatives.

An arc $i$ on $R$ is called \textit{properly embedded} if
$\partial i \subseteq \partial R$ and $i$ is transversal to
$\partial R$. $i$ is called \textit{nontrivial} (or
\textit{essential}) if $i$ cannot be deformed into $\partial R$ in
such a way that the endpoints of $i$ stay in $\partial R$ during
the deformation. The \textit{complex of arcs}, $\mathcal{B}(R)$,
on $R$ is an abstract simplicial complex. Its vertices are the
isotopy classes of nontrivial properly embedded arcs $i$ in $R$. A
set of vertices forms a simplex if these vertices can be
represented by pairwise disjoint arcs.

A nontrivial circle $a$ on a closed, connected, orientable surface
$S$ is called \textit{k-separating} (or a \textit{genus k circle})
if the surface $S_{a}$ is disconnected and one of its components
is a genus $k$ surface where $ 1 \leq k < g$. If $S_a$ is
connected, then $a$ is called \textit{nonseparating}. Two circles
on $S$ are topologically equivalent if and only if they are either
both nonseparating or $k$-separating for some $k$.

We assume that $S$ is a closed, connected, orientable surface of
genus $ g \geq 3$ throughout the paper.

\section{Superinjective Simplicial Maps of Complexes of Curves}
\begin{defn}
A simplicial map $\lambda : \mathcal{C}(S) \rightarrow
\mathcal{C}(S)$ is called {\bf superinjective} if the following
condition holds: if $\alpha, \beta$ are two vertices in
$\mathcal{C}(S)$ such that $i(\alpha,\beta) \neq 0$, then
$i(\lambda(\alpha),\lambda(\beta)) \neq 0$.
\end{defn}

\begin{lemma}
\label{injective} A superinjective simplicial map, $\lambda : \mathcal{C}(S)
\rightarrow \mathcal{C}(S)$, is injective.
\end{lemma}

\begin{proof} Let $\alpha$ and $\beta$ be two distinct vertices in $\mathcal{C}(S)$.
If $i(\alpha, \beta) \neq 0$, then $i(\lambda(\alpha),
\lambda(\beta)) \neq 0$, since $\lambda$ preserves
nondisjointness. So, $\lambda(\alpha) \neq \lambda(\beta)$. If
$i(\alpha, \beta) = 0$, we choose a vertex $\gamma$ of
$\mathcal{C}(S)$ such that $i(\gamma, \alpha)= 0$ and $i(\gamma,
\beta) \neq 0$. Then, $i(\lambda(\gamma), \lambda(\alpha)) = 0$,
$i(\lambda(\gamma), \lambda(\beta)) \neq 0$. So, $\lambda(\alpha)
\neq \lambda(\beta)$. Hence, $\lambda$ is injective. \end{proof}

\begin{lemma}
\label{conbyanedge} Let $\alpha, \beta$ be two distinct vertices of
$\mathcal{C}(S)$, and let $\lambda : \mathcal{C}(S) \rightarrow \mathcal{C}(S)$
be a superinjective simplicial map. Then, $\alpha$ and $\beta$ are
connected by an edge in $\mathcal{C} ( S )$ if and only
if $\lambda(\alpha)$ and $\lambda(\beta)$ are connected by an edge in
$\mathcal{C} (S)$.
\end{lemma}

\begin{proof} Let $\alpha, \beta$ be two distinct vertices
of $\mathcal{C}(S)$. By Lemma \ref{injective}, $\lambda$ is
injective. So, $\lambda(\alpha) \neq \lambda(\beta)$. Then we
have, $\alpha$ and $\beta$ are connected by an edge
$\Leftrightarrow$ $i(\alpha, \beta)=0$ $\Leftrightarrow$
$i(\lambda(\alpha), \lambda(\beta))=0$ (by superinjectivity)
$\Leftrightarrow$ $\lambda(\alpha)$ and $\lambda(\beta)$ are
connected by an edge.
\end{proof}

Let $P$ be a set of pairwise disjoint circles on $S$. $P$ is
called a {\it pair of pants decomposition} of $S$, if $S_P$ is a
disjoint union of genus zero surfaces with three boundary
components, pairs of pants. A pair of pants of a pants
decomposition is the image of one of these connected components
under the quotient map $q:S_P \rightarrow S$ together with the
image of the boundary components of this component. The image of
the boundary of this component is called the \textit{boundary of
the pair of pants}. A pair of pants is called \textit{embedded} if
the restriction of $q$ to the corresponding component of $S_P$ is
an embedding.

\begin{lemma}
\label{imageofpantsdecomp} Let $\lambda : \mathcal{C}(S) \rightarrow \mathcal{C}(S)$ be a
superinjective simplicial map. Let $P$ be a pair of pants decomposition of $S$. Then,
$\lambda$ maps the set of isotopy classes of elements of $P$ to the set of isotopy
classes of elements of a pair of pants decomposition, $P'$, of $S$.
\end{lemma}

\begin{proof} The set of isotopy classes of elements of $P$ forms
a top dimensional simplex,
$\bigtriangleup$, in $\mathcal{C}(S)$. Since $\lambda$ is
injective, it maps $\bigtriangleup$ to a top dimensional simplex
$\bigtriangleup'$ in $\mathcal{C}(S)$. Pairwise disjoint
representatives of the vertices of $\bigtriangleup'$ give a pair
of pants decomposition $P'$ of $S$.\end{proof}

By Euler characteristic arguments it can be seen that there exist
exactly $3g-3$ circles and $2g-2$ pairs of pants in a pair of
pants decomposition of $S$. An ordered set $(a_1, ..., a_{3g-3})$
is called an {\it ordered pair of pants decomposition} of $S$ if
$\{a_1, ..., a_{3g-3}\}$ is a pair of pants decomposition of $S$.
Let $P$ be a pair of pants decomposition of $S$. Let $a$ and $b$
be two distinct elements in $P$. Then, $a$ is called {\it
adjacent} to $b$ w.r.t. $P$ iff there exists a pair of pants in
$P$ which has $a$ and $b$ on its boundary.

{\bf Remark}: Let $P$ be a pair of pants decomposition of $S$. Let
$[P]$ be the set of isotopy classes of elements of $P$. Let
$\alpha, \beta \in [P]$. We say that $\alpha$ is adjacent to
$\beta$ w.r.t. $[P]$ if the representatives of $\alpha$ and
$\beta$ in $P$ are adjacent w.r.t. $P$. By Lemma
\ref{imageofpantsdecomp}, $\lambda$ gives a correspondence on the
isotopy classes of elements of pair of pants decompositions of
$S$. $\lambda([P])$ is the set of isotopy classes of a pair of
pants decomposition which corresponds to $P$, under this
correspondence.

\begin{lemma}
\label{adjacent} Let $\lambda : \mathcal{C}(S) \rightarrow
\mathcal{C}(S)$ be a superinjective simplicial map. Let $P$ be a
pair of pants decomposition of $S$. Then, $\lambda$ preserves the
adjacency relation for two circles, i.e. if $a, b \in P$ are two
circles which are adjacent w.r.t. $P$ and $[a]=\alpha, [b]=\beta$,
then $\lambda(\alpha), \lambda(\beta)$ are adjacent w.r.t.
$\lambda([P])$.\end{lemma}

\begin{proof} Let $P$ be a pair of pants decomposition of $S$.
Let $a, b$ be two adjacent circles in $P$ and $[a]=\alpha$,
$[b]=\beta$. Let $P_o$ be a pair of pants of $P$, having $a$ and
$b$ on its boundary. By Lemma \ref{imageofpantsdecomp}, we can
choose a pair of pants decomposition, $P'$, such that
$\lambda([P])= [P'].$

Either $P_o$ is embedded or nonembedded. In the case $P_o$ is
embedded, either $a$ and $b$ are the boundary components of
another pair of pants or not. In the case $P_o$ is nonembedded,
either $a$ or $b$ is a separating curve on $S$. So, there are four
possible cases for $P_o$. For each of these cases, in Figure 1, we
show how to choose a circle $c$ which essentially intersects $a$
and $b$ and does not intersect any other circle in $P$.

\begin{figure}[htb]
\begin{center}
\epsfxsize=1.57in \epsfbox{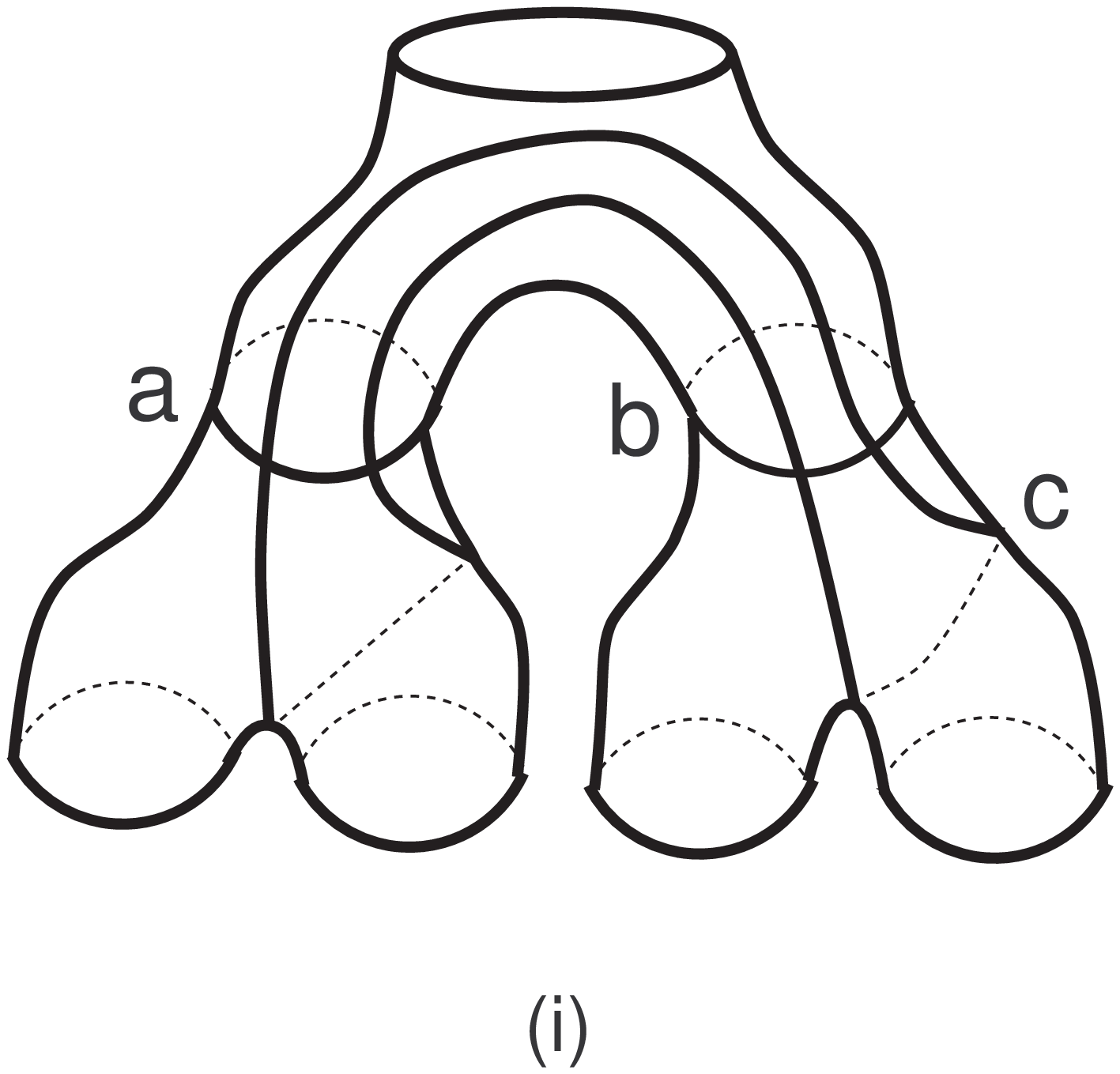} \epsfxsize=1.65in
\epsfbox{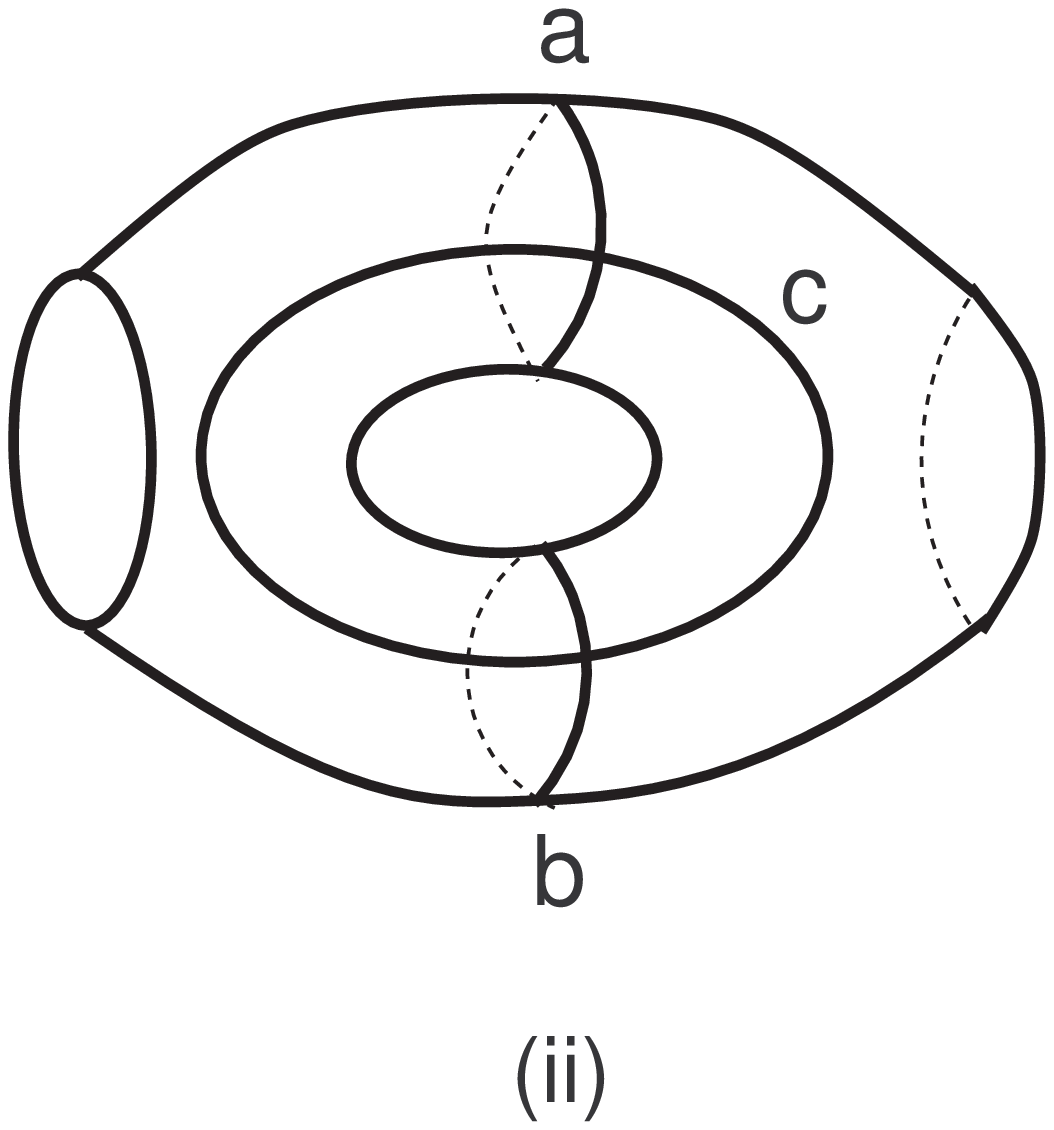}

\vspace{0.2cm} \epsfxsize=4.226in \epsfbox{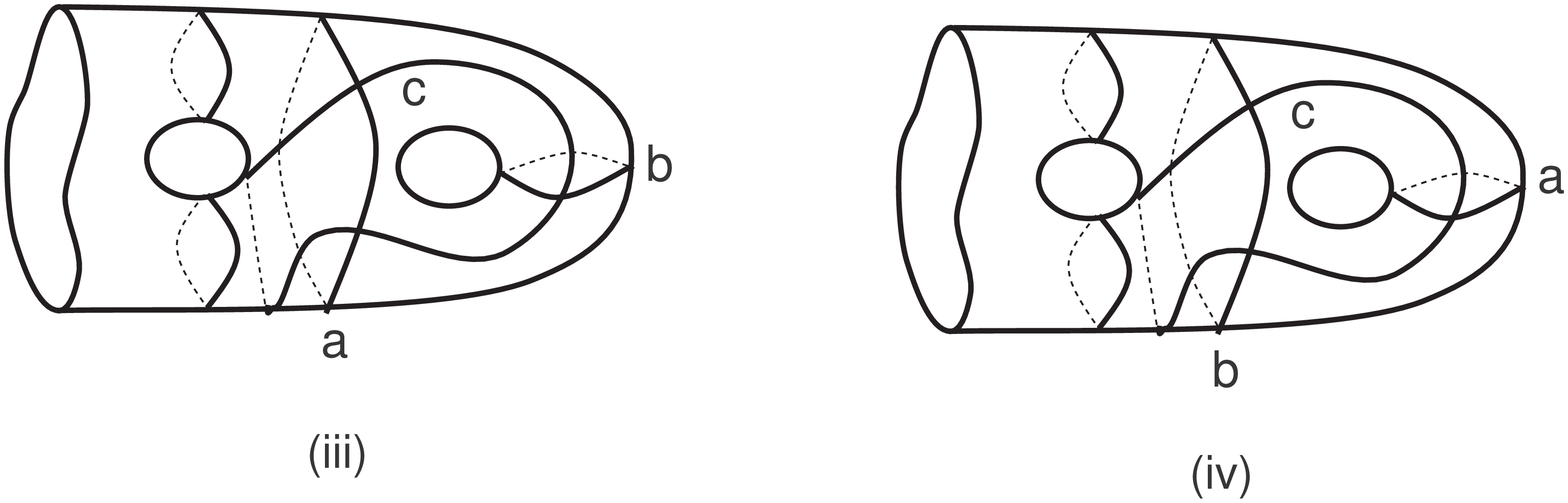}
\caption{Four possible cases for $P_o$} \label{picture1}
\end{center}
\end{figure}

Let $\gamma= [c]$. Assume that $\lambda(\alpha)$ and
$\lambda(\beta)$ do not have adjacent representatives. Since
$i(\gamma, \alpha)\neq 0 \neq i(\gamma, \beta)$, we have
$i(\lambda (\gamma), \lambda(\alpha)) \neq 0 \neq
i(\lambda(\gamma), \lambda(\beta))$ by superinjectivity. Since
$i(\gamma, [e]) = 0$ for all $e$ in $P \setminus \{a,b\}$, we have
$i(\lambda(\gamma), \lambda([e]))=0$ for all $e$ in $P \setminus
\{a,b\}$. But this is not possible because $\lambda(\gamma)$ has
to intersect geometrically essentially with some isotopy class
other than $\lambda(\alpha)$ and $\lambda(\beta)$ in the image
pair of pants decomposition to be able to make essential
intersections with $\lambda(\alpha)$ and $\lambda(\beta)$. This
gives a contradiction to the assumption that $\lambda(\alpha)$ and
$\lambda(\beta)$ do not have adjacent representatives. \end{proof}

Let $P$ be a pair of pants decomposition of $S$. A curve $x \in P$
is called a {\it 4-curve} in $P$, if there exist four distinct
circles in $P$, which are adjacent to $x$ w.r.t. $P$.

\begin{lemma}
\label{separating} Let $\lambda : \mathcal{C}(S) \rightarrow
\mathcal{C}(S)$ be a superinjective simplicial map. Then,
$\lambda$ sends the isotopy class of a $k$-separating circle to
the isotopy class of a $k$-separating circle, where $ 1 \leq k
\leq g-1$. \end{lemma}

\begin{proof}
Let $\alpha = [a]$ where $a$ is a $k$-separating circle where $1
\leq k \leq g-1$. Since the genus of $S$ is at least 3, $a$ is a
separating curve of genus at least 2. So, it is enough to consider
the cases when $k \geq 2$.

Case 1: Assume that $k = 2$. Let $S_2$ be a subsurface of $S$ of
genus $2$ having $a$ as its boundary. We can choose a pair of
pants decomposition $Q= \{a_1, a_2, a_3, a_4 \}$ of $S_2$ as shown
in Figure 2, (i). Then, we can complete $Q \cup \{a\}$ to a pair
of pants decomposition $P$ of $S$ in any way we like. By Lemma
\ref{imageofpantsdecomp}, we can choose a pair of pants
decomposition, $P'$, of $S$ such that $\lambda([P]) = [P']$.

Let $a_i'$ be the representative of $\lambda([a_i])$ which is in
$P'$ for $i=1,..,4$ and $a'$ be the representative of
$\lambda([a])$ which is in $P'$. Since $a_1$ and $a_3$ are
4-curves in $P$, by Lemma \ref{adjacent} and Lemma
\ref{injective}, $a_1'$ and $a_3'$ are 4-curves in $P'$. Notice
that a curve in a pair of pants decomposition can be adjacent to
at most four curves w.r.t. the pair of pants decomposition.

\begin{figure}[htb]
\label{picture2}
\begin{center}
\epsfxsize=1.85in \epsfbox{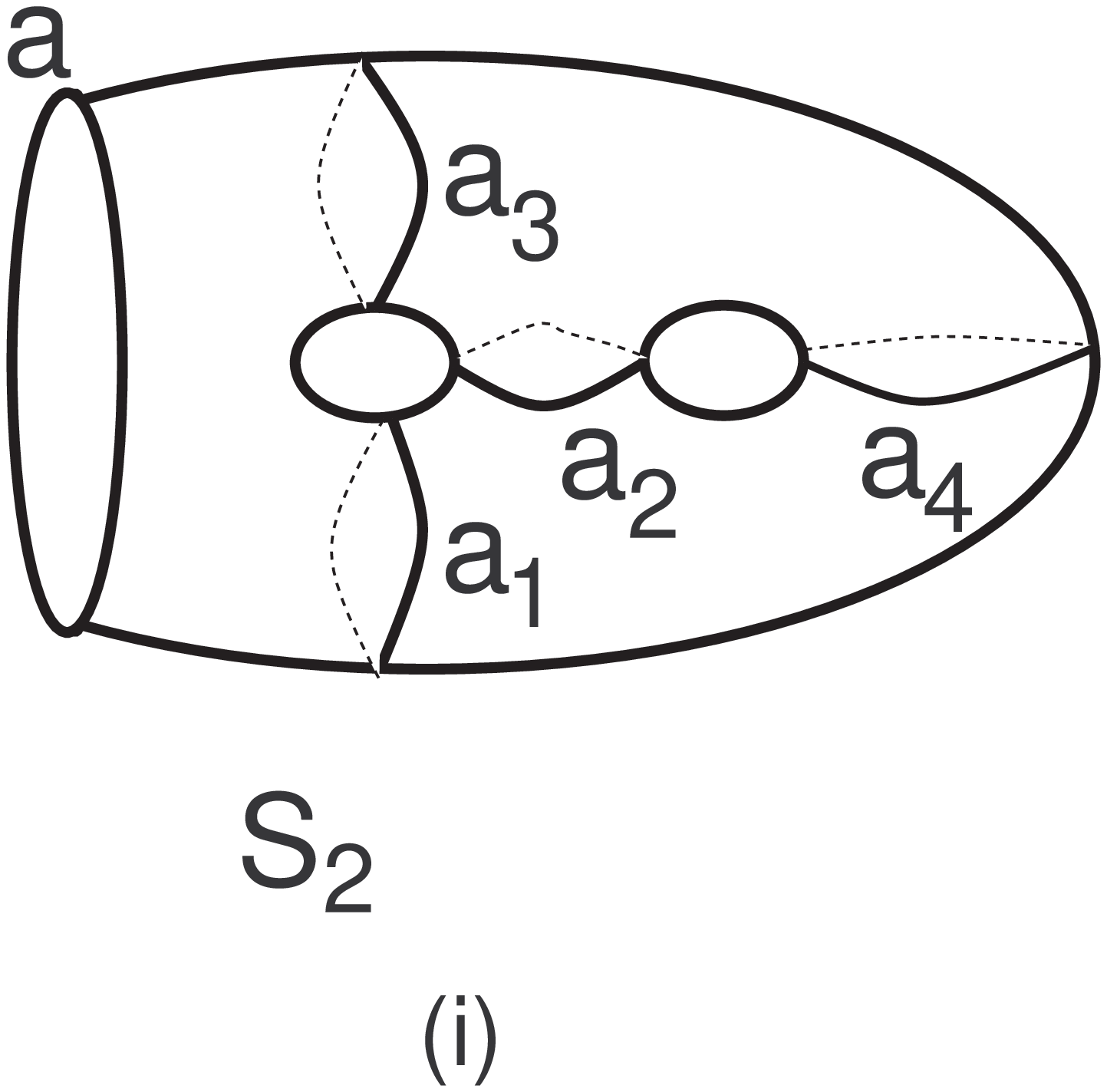} \epsfxsize=2.85in
\epsfbox{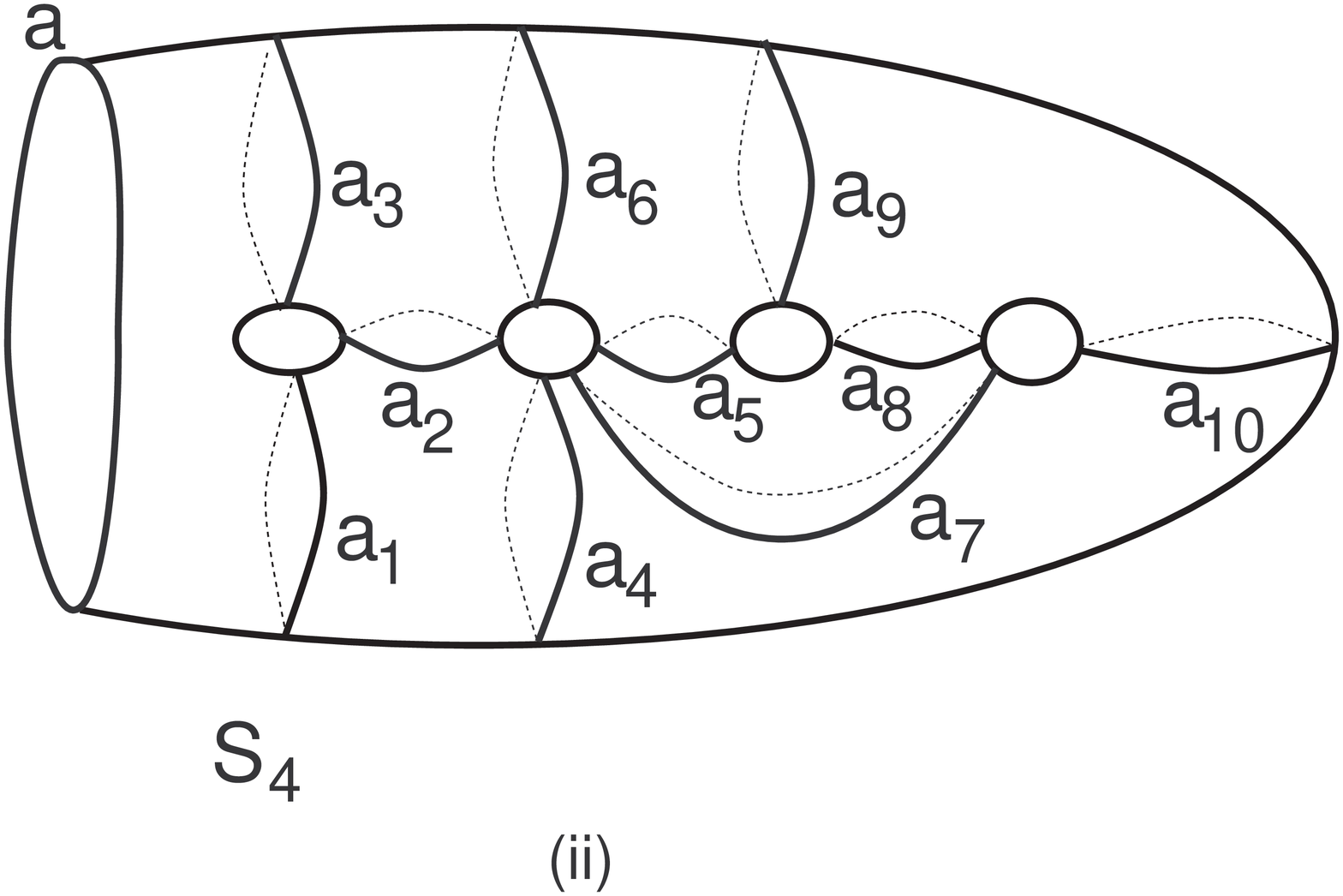} \caption{Pants decompositions with
separating circles}
\end{center}
\end{figure}

$a_3'$ is a 4-curve. So, there are two pairs of pants, $A$ and
$B$, of $P'$ such that each of them has $a_3'$ on its boundary.
Let $x, y, z, t$ be as shown in Figure 3. Since $a_3'$ is a
4-curve and $x, y, z, t$ are the only curves which are adjacent to
$a_3'$ w.r.t. $P'$, $x, y, z, t$ are four distinct curves.

$a_3$ is adjacent to $a_1$ w.r.t. $P$ implies that $a_3'$ is
adjacent to $a_1'$ w.r.t. $P'$. Then, W.L.O.G. we can assume that
$x = a_1'$. So, $a_1'$ is a boundary component of $A$. Since
$a_1'$ is a 4-curve in $P'$, $a_1'$ is also a boundary component
of a pair of pants different from $A$. Since $a_1' \neq z$, $a_1'
\neq t$ and $a_1' \neq a_3'$, $a_1'$ is not a boundary component
of $B$. So, there is a new pair of pants, $C$, which has $a_1'$ on
its boundary. Let $v, w$ be as shown in Figure 3. Since $a_1'$ is
a 4-curve and $y, v, w, a_3'$ are the only curves which are
adjacent to $a_1'$ w.r.t. $P'$, $y, v, w, a_3'$ are four distinct
curves. Since $a_1$ is adjacent to each of $a, a_2, a_3, a_4$
w.r.t. $P$, $a_1'$ is adjacent to each of $a', a_2', a_3', a_4'$
w.r.t. $P'$. Then, $\{y, v, w, a_3'\} = \{a', a_2', a_3', a_4'\}$
and, so, $\{y, v, w\} = \{a', a_2', a_4'\}$. Similarly, since
$a_3$ is adjacent to each of $a, a_1, a_2, a_4$ w.r.t. $P$, $a_3'$
is adjacent to each of $a', a_1', a_2', a_4'$ w.r.t. $P'$. Then,
since each of $x, y, z, t$ is adjacent to $a_3'$, we have $\{x, y,
z, t\} = \{a', a_1', a_2', a_4'\}$. Since $x= a_1'$, $\{y, z, t\}
= \{a', a_2', a_4'\} = \{y, v, w\}$. Then $\{z, t\} = \{v, w\}$.
Hence, $A \cup B \cup C$ is a genus 2 subsurface of $S$ having $y$
as its boundary.

\begin{figure}[htb]
\label{picture3}
\begin{center}
\epsfig{file=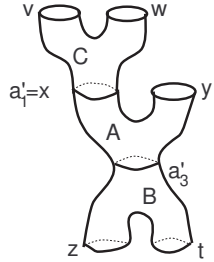,width=3.9cm} \caption{Adjacent
circles}
\end{center}
\end{figure}

Since $\{y, v, w\} = \{a', a_2', a_4'\}$, $y$ is $a'$ or $a_2'$ or
$a_4'$. Assume that $y=a_2'$. Then, $a_2'$ can't be adjacent to
$a_4'$ w.r.t. $P'$ as it can be seen from Figure 3. But this gives
a contradiction since $a_2$ is adjacent to $a_4$ w.r.t. $P$,
$a_2'$ must also be adjacent to $a_4'$ w.r.t. $P'$. Similarly, if
$y = a_4'$, we get a contradiction since $a_4'$ wouldn't be
adjacent to $a_2'$ w.r.t. $P'$, but it should be as $a_4$ is
adjacent to $a_2$ w.r.t. $P$. Hence, $y=a'$. This proves that $a'$
is a genus 2 curve.

Case 2: Assume that $k \geq 3$. Let $S_k$ be a subsurface of $S$
of genus $k$ having $a$ as its boundary. We can choose a pair of
pants decomposition $Q = \{a_1, a_2, ..., a_{3k-2} \}$ of $S_k$,
and then complete $Q \cup \{a\}$ to a pair of pants decomposition
$P$ of $S$ such that each of $a_i$ is a 4-curve in $P$ for
$i=1,2,...,3k-2$ and $a, a_1, a_3$ are the boundary components of
a pair of pants of $Q$. In Figure 2, (ii), we show how to choose
$Q$ when $k=4$. In the other cases, when $k=3$ or $k \geq 5$, a
similar pair of pants decomposition of $S_k$ can be chosen.

Let $P'$ be a pair of pants decomposition of $S$ such that
$\lambda([P]) = [P']$. Let $a_i'$ be the representative of
$\lambda([a_i])$ which is in $P'$, for $i=1,..., 3k-2$, and $a'$
be the representative of $\lambda([a])$ which is in $P'$. Let $Q'=
\{a_1', a_2', ..., a_{3k-2}'\}.$

For every 4-curve $x$ in $P'$, there exist two pairs of pants
$A(x)$ and $B(x)$ of $P'$ having $x$ as one of their boundary
components. Let $C(x)= A(x) \cup B(x)$. The boundary of $C(x)$
consists of four distinct curves, which are adjacent to $x$ w.r.t.
$P'$. Since all the curves $a_1, ..., a_{3k-2}$ in $Q$ are
4-curves in $P$, all the corresponding curves $a_1',...,
a_{3k-2}'$ are 4-curves in $P'$ by Lemma \ref{adjacent} and Lemma
\ref{injective}. Let $N'= C(a_1') \cup C(a_2') \cup ... \cup
C(a_{3k-2}')$.

Claim: $N'$ is a genus $k$ subsurface of $S$ having $a'$ as its
boundary.

Proof: Each boundary component of $C(a_i')$ is either $a'$ or
$a_j'$ for some $j=1, ..., 3k-2$. Since $a_i'$ is in the interior
of $C(a_i')$, $a_i'$ is in the interior of $N'$ for
$i=1,2,...,3k-2$. So, all the boundary components of $C(a_i')$
which are different from $a'$ are in the interior of $N'$. Hence,
$N'$ has at most one boundary component, which could be $a'$. We
know that $a'$ is adjacent to two distinct curves, $a_1', a_3'$,
w.r.t. $P'$ since $a$ is adjacent to $a_1, a_3$ w.r.t. $P$.
Suppose that $a'$ is in the interior of $N'$. Then $a'$ is only
adjacent to curves w.r.t. $P'$ which are in $Q'$. On the other
hand, being adjacent to $a_1', a_3'$, it has to be adjacent to a
third curve, $a_j'$ in $Q'$, w.r.t. $P'$. But each of such curves
has four distinct adjacent curves w.r.t. $P'$ which are different
from $a'$ already since for $j \neq 1, 3$, each of $a_j$ has four
distinct adjacent curves w.r.t. $P$, which are different from $a$
by our choice of $P$. So, $a'$ cannot be adjacent to $a_j'$ w.r.t.
$P'$ when $j \neq 1,3$. Hence, $a'$ is not in the interior of
$N'$. It is on the boundary of $N'$. So, $N'$ is a subsurface of
$S$ having $a'$ as its boundary. To see that $N'$ has genus $k$,
it is enough to realize that $\{a_1',..., a_{3k-2}'\}$ is a pair
of pants decomposition of $N'$. Hence, $a'$ is a $k$-separating
circle.\end{proof}

\begin{lemma}
\label{division} Let $\lambda : \mathcal{C}(S) \rightarrow
\mathcal{C}(S)$ be a superinjective simplicial map. Let $t$ be a
$k$-separating circle on $S$, where $1 \leq k \leq g-1$. Let $S_1,
S_2$ be the distinct subsurfaces of $S$ of genus $k$ and $g-k$
respectively which have $t$ as their boundary. Let $t' \in
\lambda([t])$. Then $t'$ is a $k$-separating circle and there
exist subsurfaces $S_1', S_2'$ of $S$ of genus $k$ and $g-k$
respectively which have $t'$ as their boundary such that
$\lambda(\mathcal{C}(S_1)) \subseteq \mathcal{C}(S_1')$ and
$\lambda(\mathcal{C}(S_2)) \subseteq \mathcal{C}(S_2')$.
\end{lemma}

\begin{proof}
Let $t$ be a $k$-separating circle where $1 \leq k \leq g-1$.
Since the genus of $S$ is at least 3, $t$ is a separating curve of
genus at least 2. So, it is enough to consider the cases when $k
\geq 2$.

Let $S_1, S_2$ be the distinct subsurfaces of $S$ of genus $k$ and
$g-k$ respectively which have $t$ as their boundary. Let $t' \in
\lambda([t])$. By Lemma \ref{separating}, $t'$ is a $k$-separating
circle. As we showed in the proof of Lemma \ref{separating}, there
is a pair of pants decomposition $P_1$ of $S_1$, and $P_1 \cup
\{t\}$ can be completed to a pair of pants decomposition $P$ of
$S$ such that a set of curves, $P_1'$, corresponding (via
$\lambda$) to the curves in $P_1$, can be chosen such that $P_1'$
is a pair of pants decomposition of a subsurface that has $t'$ as
its boundary. Let $S_1'$ be this subsurface. Let $S_2'$ be the
other subsurface of $S$ which has $t'$ as its boundary. A pairwise
disjoint representative set, $P'$, of $\lambda([P])$ containing
$P_1' \cup \{t'\}$ can be chosen. Then, by Lemma
\ref{imageofpantsdecomp}, $P'$ is a pair of pants decomposition of
$S$. Let $P_2 = P \setminus (P_1 \cup {t})$ and $P_2' = P'
\setminus (P_1' \cup {t'})$. Then $P_2, P_2'$ are pair of pants
decompositions of $S_2, S_2'$ respectively as $P_1, P_1'$ are pair
of pants decompositions of $S_1, S_1'$ respectively.

Now, let $\alpha$ be a vertex in $\mathcal{C}(S_1)$. Then, either
$\alpha \in [P_1]$ or $\alpha$ has a nonzero geometric
intersection with an element of $[P_1]$. In the first case,
clearly $\lambda(\alpha) \in \mathcal{C}(S_1')$ since elements of
$[P_1]$ correspond to elements of $[P_1'] \subseteq
\mathcal{C}(S_1')$. In the second case, since $\lambda$ preserves
zero and nonzero geometric intersection (since $\lambda$ is
superinjective) and $\alpha$ has zero geometric intersection with
the elements of $[P_2]$ and $[t]$ and nonzero intersection with an
element of $[P_1]$, $\lambda(\alpha)$ has zero geometric
intersection with elements of $[P_2']$ and $[t']$ and nonzero
intersection with an element of $[P_1']$. Then, $\lambda(\alpha)
\in \mathcal{C}(S_1')$. Hence, $\lambda(\mathcal{C}(S_1))
\subseteq \mathcal{C}(S_1')$. The proof of
$\lambda(\mathcal{C}(S_2)) \subseteq \mathcal{C}(S_2')$ is
similar.\end{proof}

\begin{lemma}
\label{top} Let $\lambda : \mathcal{C}(S) \rightarrow
\mathcal{C}(S)$ be a superinjective simplicial map. Then $\lambda$
preserves topological equivalence of ordered pairs of pants
decompositions of $S$, (i.e. for a given ordered pair of pants
decomposition $P=(c_1, c_2, ..., c_{3g-3})$ of $S$, and a
corresponding ordered pair of pants decomposition $P'=(c_1', c_2',
..., c_{3g-3}')$ of $S$, where $[c_i']= \lambda([c_i])$ $\forall
i= 1, 2, ..., 3g-3$, there exists a homeomorphism $F: S
\rightarrow S$ such that $F(c_i)=c_i'$ $\forall i= 1, 2, ...,
3g-3$).\end{lemma}

\begin{proof} Let $P$ be a pair of pants decomposition of $S$ and
$A$ be a nonembedded pair of pants in $P$. The boundary of $A$
consists of the circles $x, y$ where $x$ is a 1-separating circle
on $S$ and $y$ is a nonseparating circle on $S$. Let $R$ be the
subsurface of $S$ of genus $g-1$ which is bounded by $x$ and let
$T$ be the subsurface of $S$ of genus 1 which is bounded by $x$.
Let $P_1$ be the set of elements of $P \setminus \{x\}$ which are
on $R$ and $P_2$ be the set of elements of $P \setminus \{x\}$
which are on $T$. Then, $P_1, P_2$ are pair of pants
decompositions of $R, T$ respectively. So, $P_2 = \{y\}$ is a pair
of pants decomposition of $T$. By Lemma \ref{division}, there
exists a 1-separating circle $x' \in \lambda([x])$ and subsurfaces
$T', R'$, of $S$, of genus 1 and $g-1$ respectively such that
$\lambda(\mathcal{C}(R)) \subseteq \mathcal{C}(R')$ and
$\lambda(\mathcal{C}(T)) \subseteq \mathcal{C}(T')$. Since $[P_1]
\subseteq \mathcal{C}(R)$, we have $\lambda([P_1]) \subseteq
\mathcal{C}(R')$. Since $[P_2] \subseteq \mathcal{C}(T)$, we have
$\lambda([P_2]) \subseteq \mathcal{C}(T')$. Since $\lambda$
preserves disjointness, we can see that a set, $P_1'$, of pairwise
disjoint representatives of $\lambda([P_1])$ disjoint from $x'$
can be chosen. By counting the number of curves in $P_1'$, we can
see that $P_1'$ is a pair of pants decomposition of $R'$.
Similarly, a set, $P_2'$, of pairwise disjoint representatives of
$\lambda([P_2])$ disjoint from $x'$ can be chosen. By counting the
number of curves in $P_2'$,  we can see that $P_2'$ is a pair of
pants decomposition of $T'$. Since $P_2$ has one element, $y$,
$P_2'$ has one element. Let $y' \in P_2'$. Since $x', y'$
correspond to $x,y$ respectively and $y$ and $y'$ give pair of
pants decompositions on $T$ and $T'$ (which are both nonembedded
pairs of pants) and $x$ and $x'$ are the boundaries of $R$ and
$R'$, we see that $\lambda$ ``sends" a nonembedded pair of pants
to a nonembedded pair of pants.

Let $B$ be an embedded pair of pants of $P$. Let $x, y, z \in P$
be the boundary components of $B$. We consider two cases:

(i) At least one of $x$, $y$ or $z$ is a separating circle.\\
\indent (ii) All of $x, y, z$ are nonseparating circles.

In the first case, W.L.O.G assume that $x$ is a $k$-separating
circle for $1 \leq k < g$. Let $S_1, S_2$ be the distinct
subsurfaces of $S$ of genus $k$ and $g-k$ respectively which have
$x$ as their boundary. W.L.O.G. assume that $y, z$ are on $S_2$.
Let $x' \in \lambda([x])$. By Lemma \ref{division}, there exist
subsurfaces, $S_1', S_2'$, of $S$ of genus $k$ and $g-k$
respectively which have $x'$ as their boundary such that
$\lambda(\mathcal{C}(S_1)) \subseteq \mathcal{C}(S_1')$ and
$\lambda(\mathcal{C}(S_2)) \subseteq \mathcal{C}(S_2')$. Then,
since $y \cup z \subseteq S_2$, $\lambda(\{[y], [z]\}) \subseteq
\mathcal{C}(S_2')$. Let $y' \in \lambda([y]), z' \in \lambda([z])$
such that $\{x', y', z'\}$ is pairwise disjoint. Let $P'$ be a set
of pairwise disjoint representatives of $\lambda([P])$ which
contains $x', y', z'$. $P'$ is a pair of pants decomposition of
$S$. Then, since $x$ is adjacent to $y$ and $z$ w.r.t. $P$, $x'$
is adjacent to $y'$ and $z'$ w.r.t. $P'$ by Lemma \ref{adjacent}.
Then, since $x' \cup y'\cup z' \subseteq S_2'$, and $x'$ is the
boundary of $S_2'$, there is an embedded pair of pants in $S_2'$
which has $x', y', z'$ on its boundary. So, $\lambda$ ``sends'' an
embedded pair of pants bounded by $x, y, z$ to an embedded pair of
pants bounded by $x', y', z'$ in this case.

In the second case, we can find a nonseparating circle $w$ and a
2-separating circle $t$ on $S$ such that $\{x, y, z, w\}$ is
pairwise disjoint and $x, y, z, w$ are on a genus 2 subsurface,
$S_1$, that $t$ bounds as shown in Figure 4. Let $P_1 = \{x, y, z,
w \}$. $P_1$ is a pair of pants decomposition of $S_1$. We can
complete $P_1 \cup \{t\}$  to a pants decomposition $P$ of $S$.

\begin{figure}[htb]
\label{picture4}

\begin{center}
\epsfig{file=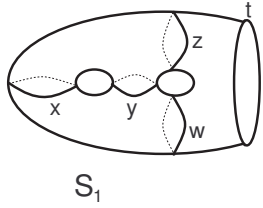,width=4.5cm} \caption{Nonseparating
circles, $x, y, z$, bounding a pair of pants}
\end{center}
\end{figure}

Let $S_2$ be the subsurface of $S$ of genus $g-2$ which is not
equal to $S_1$ and that is bounded by $t$. By Lemma
\ref{division}, there exist a 2-separating circle $t' \in
\lambda([t])$ and subsurfaces, $S_1', S_2'$, of $S$ of genus 2 and
$g-2$ respectively such that $\lambda(\mathcal{C}(S_1)) \subseteq
\mathcal{C}(S_1')$ and $\lambda(\mathcal{C}(S_2)) \subseteq
\mathcal{C}(S_2')$. Since $P_1 \subseteq S_1$, $\lambda([P_1])
\subseteq \mathcal{C}(S_1')$. We can choose a set, $P_1'$, of
pairwise disjoint representatives of $\lambda([P_1])$ on $S_1'$.
Then, $P_1' \cup \{t'\}$ is a pair of pants decomposition of
$S_1'$. We can choose a pairwise disjoint representative set,
$P'$, of $\lambda([P])$ containing $P_1'$. $P'$ is a pair of pants
decomposition of $S$. Let $x', y', z', w' \in P_1'$ be the
representatives of $x, y, z, w$ respectively. Then, since $t$ is
adjacent to $z$ and $w$ w.r.t. $P$, $t'$ is adjacent to $z'$ and
$w'$ w.r.t. $P'$ by Lemma \ref{adjacent}. Then, since $t' \cup z'
\cup w' \subseteq S_1'$ and $t'$ is the boundary of $S_1'$, there
is an embedded pair of pants in $S_1'$ which has $t', z', w'$ on
its boundary. Since $z$ is a 4-curve in $P$, $z'$ is a 4-curve in
$P'$. Since $z$ is adjacent to $x, y$ w.r.t. $P$, $z'$ is adjacent
to $x', y'$ w.r.t. $P'$. Since $z'$ is on the boundary of a pair
of pants which has $w', t'$ on its boundary, and $z'$ is adjacent
to $x', y'$, there is a pair of pants having $x', y', z'$ on its
boundary. So, $\lambda$ ``sends'' an embedded pair of pants
bounded by $x, y, z$ to an embedded pair of pants bounded by $x',
y', z'$ in this case too.

Assume that $P = (c_1, c_2, ..., c_{3g-3})$ is an ordered pair of
pants decomposition of $S$. Let $c_i' \in \lambda([c_i])$ such
that the elements of $\{c_1', c_2', ..., c_{3g-3}'\}$ are pairwise
disjoint. Then, $P'=(c_1', c_2', ..., c_{3g-3}')$ is an ordered
pair of pants decomposition of $S$. Let $(B_1, B_2, ...,
B_{2g-2})$ be an ordered set containing the connected components
of $S_P$. By the arguments given above, there is a corresponding,
``image'', ordered collection of pairs of pants $(B_1', B_2',...,
B_{2g-2}')$. Nonembedded pairs of pants correspond to nonembedded
pairs of pants and embedded pairs of pants correspond to embedded
pairs of pants. Then, by the classification of surfaces, there
exists an orientation preserving homeomorphism $h_i : B_i
\rightarrow B_i'$, for all $i =1,..., 2g-2$. We can compose each
$h_i$ with an orientation preserving homeomorphism $r_i$ which
switches the boundary components, if necessary, to get $h_i'= r_i
\circ h_i$ to agree with the correspondence given by $\lambda$ on
the boundary components, (i.e. for each boundary component $a$ of
$B_i$ for $i=1,...,2g-2$, $\lambda([q(a)])= [q'(h_i'(a))]$ where
$q: S_P \rightarrow S$ and $q': S_{P'} \rightarrow S$ are the
natural quotient maps). Then for two pairs of pants with a common
boundary component, we can glue the homeomorphisms by isotoping
the homeomorphism of the one pair of pants so that it agrees with
the homeomorphism of the other pair of pants on the common
boundary component. By adjusting these homeomorphisms on the
boundary components and gluing them we get a homeomorphism $F : S
\rightarrow S$ such that $F(c_i)=c_i'$ for all $i=1,2, ...,
3g-3$.\end{proof}

{\bf Remark:} Let $\mathcal{E}$ be an ordered set of vertices of
$\mathcal{C}(S)$ having a pairwise disjoint representative set
$E$. Then, $E$ can be completed to an ordered pair of pants
decomposition, $P$, of $S$. We can choose an ordered pairwise
disjoint representative set, $P'$, of $\lambda([P])$ by Lemma
\ref{imageofpantsdecomp}. Let $E'$ be the elements of $P'$ which
correspond to the elements of $E$. By Lemma \ref{top}, $P$ and
$P'$ are topologically equivalent as ordered pants decompositions.
Hence, the set $E$ and $E'$ are topologically equivalent. So,
$\lambda$ gives a correspondence which preserves topological
equivalence on a set which has pairwise disjoint representatives.
In particular, $\lambda$ ``sends'' the isotopy class of a
nonseparating circle to the isotopy class of a nonseparating
circle.

We use the following lemma to understand some more properties of
superinjective simplicial maps.

\begin{lemma} {\bf (Ivanov) \cite{Iv1}}
\label{Ivanovlemma} Let $\alpha_{1}$ and $\alpha_{2}$ be two
vertices in $\mathcal{C}(S)$. Then, $i( \alpha_{1}, \alpha_{2})=1$
if and only if there exist isotopy classes $\alpha_{3},
\alpha_{4}, \alpha_{5}$ such that\\
\indent (i) $i(\alpha_{i}, \alpha_{j})=0$ if and only if the
$i^{th}$
and $j^{th}$ circles on Figure 5 are disjoint. \\
\indent (ii) if $\alpha_{4}$ is the isotopy class of a circle
$C_{4}$, then $C_4$ divides $S$ into two pieces, and one of these
is a torus with one hole containing some representatives of the
isotopy classes $\alpha_1, \alpha_2$. \end{lemma}

\begin{figure}
\begin{center}
\epsfxsize=2.3in \epsfbox{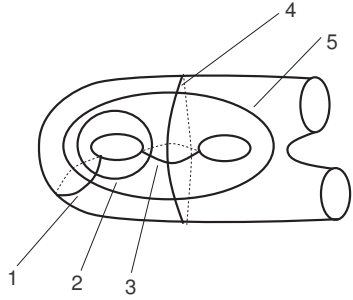} \caption{Circles
intersecting once}
\end{center}
\end{figure}

\begin{lemma}
\label{intone} Let $\lambda : \mathcal{C}(S) \rightarrow \mathcal{C}(S)$ be a superinjective
simplicial map. Let $\alpha$, $\beta$ be two vertices of $\mathcal{C}(S)$. If $i(\alpha, \beta)=1$,
then $i(\lambda(\alpha), \lambda(\beta))=1$.
\end{lemma}

\begin{proof} Let $\alpha$, $\beta$ be two vertices of $\mathcal{C}(S)$
such that $i(\alpha, \beta)=1$. Then, by Ivanov's Lemma, there
exist isotopy classes $\alpha_{3}, \alpha_{4}, \alpha_{5}$ such
that $i(\alpha_{i}, \alpha_{j})=0$ if and only if $i^{th}, j^{th}$
circles on Figure 5 are disjoint and if $\alpha_{4}$ is the
isotopy class of a circle $C_{4}$, then $C_4$ divides $S$ into two
pieces, and one of these is a torus with one hole containing some
representative of the isotopy classes $\alpha_1, \alpha_2$. Then,
since $\lambda$ is superinjective $i(\lambda(\alpha_{i}),
\lambda(\alpha_{j}))=0$ if and only if $i^{th}, j^{th}$ circles on
Figure 5 are disjoint, and by Lemma \ref{division}, if
$\lambda(\alpha_{4})$ is the isotopy class of a circle $C_{4}'$,
then $C_4'$ divides $S$ into two pieces, and one of these is a
torus with one hole containing some representative of the isotopy
classes $\lambda(\alpha_1), \lambda(\alpha_2)$. Then, by Ivanov's
Lemma, $i(\lambda(\alpha), \lambda(\beta))=1$.
\end{proof}

\section{Induced Map On Complex Of Arcs}

In this section, we assume that $\lambda : \mathcal{C}(S)
\rightarrow \mathcal{C}(S)$ is a superinjective simplicial map,
$c, d$ are nonseparating circles on $S$, and $[d]=\lambda([c])$.
Let $\mathcal{V}(S_c)$ and $\mathcal{V}(S_d)$ be the sets of
vertices of $\mathcal{B}(S_c)$ and $\mathcal{B}(S_d)$
respectively. We prove that $\lambda$ induces a map $\lambda_* :
\mathcal{V}(S_c) \rightarrow \mathcal{V}(S_d)$ with certain
properties. Then we prove that $\lambda_*$ extends to an injective
simplicial map $\lambda_* : \mathcal{B}(S_c) \rightarrow
\mathcal{B}(S_d)$.

First, we prove some lemmas which we use to see some properties of $\lambda$.

\begin{lemma}
\label{A} Let $a$ and  $b$ be two disjoint arcs on $S_c$
connecting the two boundary components, $\partial_1, \partial_2$,
of $S_c$. Let $N$ be a regular neighborhood of $a \cup b \cup
\partial_1 \cup \partial_2$ in $S_c$. Then, $(N, a, b) \cong (S_4
^2, a_o, b_o)$ where $S_4^2$ is a standard sphere with four holes
and $a_o, b_o$ are arcs as shown in Figure 6.

\begin{figure}
\begin{center}
\epsfxsize=2.3in \epsfbox{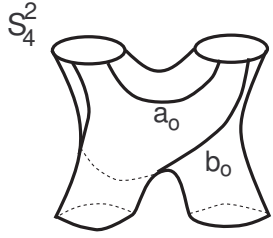} \caption{Unlinked
disjoint arcs of type 2}
\end{center}
\end{figure}
\end{lemma}

\begin{proof}  Let $\Gamma = \partial_1 \cup \partial_2 \cup a \cup b$.
Let $N$ be a regular neighborhood of $\Gamma$ in $S_c$. $N$
deformation retracts onto $\Gamma$. So, $N$ and $\Gamma$ have the
same Euler characteristics. Let $m$ be the genus of $N$ and $n$ be
the number of boundary components of $N$. It is easy to see that
$\chi(\Gamma)= -2$. Then, $-2=\chi(N)=2-2m-n$. So, $2m+n=4$.

\begin{figure}[htb]
\begin{center}
\hspace{0.052in} \epsfxsize=1.6in \epsfbox{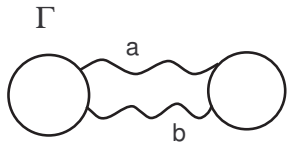}
\caption{The graph of $a, b$ and the boundary components}

\vspace{0.15in}
\end{center}
\end{figure}

There are three possibilities for $m$ and $n$.

(i) $m=0, n=4$, \hspace{0.5cm} (ii) $m=1, n=2$, \hspace{0.5cm}
(iii) $m=2, n=0$.

Since there are 2 boundary components corresponding to
$\partial_1, \partial_2$, none of which is the boundary component
on a given side of the arc $a$, $n$ is at least 3. So, only (i)
can hold. Therefore, $N$ is homeomorphic to $S_4^2$, a sphere with
four holes.

By using Euler characteristic arguments we can see that a regular
neighborhood of $\partial_1 \cup \partial_2 \cup a$ in $N$ is
homeomorphic to a pair of pants, $P$. Let $w$ be the boundary
component of $P$ different from $\partial_1, \partial_2$.

Let $\Gamma'$ be the graph that we get when we contract
$\partial_1$ and $\partial_2$ to two points. If $a$ is isotopic to
$b$ in $N$, the two arcs corresponding to $a$ and $b$ on $\Gamma'$
should be homotopic to each other relative to these two points.
But this gives a contradiction since the arcs intersect only at
these end points and the union of the arcs is $\gamma'$ which is a
circle and two such arcs cannot be homotopic relative to their end
points on a circle. So, $a$ is not isotopic to $b$ in $N$. Since
$b$ connects the two boundary components $\partial_1$ and
$\partial_2$ and $\partial_1, \partial_2 \in P$, $b \cap P$ is
nonempty. W.L.O.G. assume that $b$ intersects the boundary
components of $P$ transversely and doesn't intersect $a$. Since
$N$ is a regular neighborhood of $a \cup b \cup \partial_1 \cup
\partial_2$ and $b$ is a properly embedded essential arc which is
nonisotopic to $a$ in $N$, $b \cap P$ contains exactly two
essential properly embedded arcs, let's call them $b_1, b_2$. One
of them starts on $\partial_1$ and ends on $w$ and the other one
starts on $\partial_2$ and ends on $w$. Let $P'$ be a regular
neighborhood of $\partial_1 \cup
\partial_2 \cup a$ in $P$ such that $P' \cap b = b_1 \cup b_2$. Let $x$
be the boundary component of $P'$, which is different from
$\partial_1, \partial_2$. We choose this neighborhood to get rid
of the possible inessential arcs of $b$ in $P$. $P'$ is a pair of
pants. Since $N$ is a sphere with four holes, the complement of
$P'$ in $N$ is a pair of pants, $R$. Let $y$ and $z$ be the
boundary components of $R$ which are different from $x$. Then, we
have a homeomorphism $\phi:(N,
\partial_1, \partial_2, x, y, z) \rightarrow (S_4^2, r_o, t_o,
x_o, y_o, z_o)$ where $r_o, t_o, x_o, y_o, z_o$ are as shown in
Figure 8.

Let $P_o$ be the pair of pants bounded by $r_o, t_o, x_o$. Let
$x_1, x_2$ be two parallel copies of $x_o$ in $P_o$ as shown in
Figure 8 such that each of them intersects $\phi(b)$ transversely
at exactly 2 points and none of $x_1, x_2$ intersects $\phi(a)$.
Let $A_o$ be the annulus which is bounded by $x_1$ and $x_2$, and
$B_o$ be the annulus which is bounded by $x_o, x_1$. Let $Q_o$ be
the pair of pants bounded by $r_o, t_o, x_2$.

\begin{figure}[htb]
\begin{center}
\epsfxsize=3.4in \epsfbox{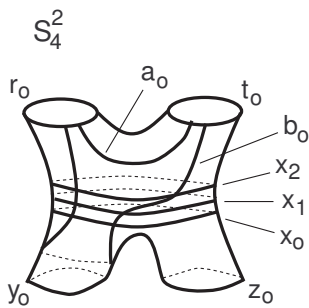}

\caption{Disjoint arcs on pairs of pants}
%Figure 8
\end{center}
\end{figure}

By the classification of properly embedded essential arcs on a
pair of pants, there exists an isotopy on $Q_o$ carrying $\phi(a)$
and $Q_o \cap \phi(b)$ to $Q_o \cap a_o$ and $Q_o \cap b_o$
respectively, where $a_o, b_o$ are as shown in the figure. Let
$\kappa : Q_o \times I \rightarrow Q_o$ be such an isotopy. We can
extend $\kappa$ to $\tilde{\kappa}: (Q_o \cup A_o) \times I
\rightarrow (Q_o \cup A_o)$ so that $\tilde{\kappa }_t$ is $id$ on
$x_1$ for all $t \in I$.

Let $R_o$ be the pair of pants bounded by $x_o, y_o, z_o$. Suppose
that $R_o \cap \phi(b)$ is an inessential arc in $R_o$. Then it
can be deformed into the interior of $P_o$ and then we get
$\phi(a)$ is isotopic to $\phi(b)$ in $P_o$ which implies that $a$
is isotopic to $b$ in $N$. This gives a contradiction. So, $b \cap
R_o$ is an essential arc in $R_o$. Then, by the classification of
properly embedded essential arcs on a pair of pants there exists
an isotopy carrying $R_o \cap \phi(b)$ to $R_o \cap b_o$. Let
$\tau : R_o \times I \rightarrow R_o$ be such an isotopy. We can
extend $\tau$ to $\tilde{\tau}:(R_o \cup B_o) \times I \rightarrow
(R_o \cup B_o)$ so that $\tilde{\tau }_t$ is $id$ on $x_1$ for all
$t \in I$. Then, by gluing the extensions $\tilde{\kappa}$ and
$\tilde{\tau}$ we get an isotopy $\vartheta$ on $S_4^2 \times I$
which fixes each of $r_o, t_o, x_o, y_o, z_o$. By the
classification of isotopy classes of arcs relative to the boundary
on an annulus, $\vartheta_1(\phi(b)) \cap (A_o \cup B_o)$ can be
isotoped to $t_{x_o} ^k(b_o) \cap (A_o \cup B_o)$ for some $k \in
\mathbb{Z}$. Let's call this isotopy $\mu$. Let $\tilde{\mu }$
denote the extension by id to $N$. Then we have, $t_{x_o}
^{-k}(\tilde{\mu}_1(\vartheta_1 (\phi(b)))) = b_o$. Clearly,
$t_{x_o} ^{-k} \circ \tilde{\mu}_1 \circ \vartheta _1$ fixes each
of $r_o, t_o, x_o, y_o, z_o$. Hence, we get a homeomorphism,
$t_{x_o}^{-k} \circ \tilde{\mu}_1 \circ \vartheta_1 \circ \phi:
(N, a, b) \rightarrow (S_4^2, a_o, b_o)$.\end{proof}

Let $a$ and $b$ be two disjoint arcs connecting a boundary
component of $S_c$ to itself. Then, $a$ and $b$ are called {\it
linked} if their end points alternate on the boundary component.
Otherwise, they are called {\it unlinked}.

The proofs of Lemma \ref{B}, Lemma \ref{C} and Lemma \ref{D} are
similar to the proof of Lemma \ref{A}. So, we do not prove these
lemmas here. We only state them.

\begin{figure}[htb]
\begin{center} \epsfxsize=5.6in
\epsfbox{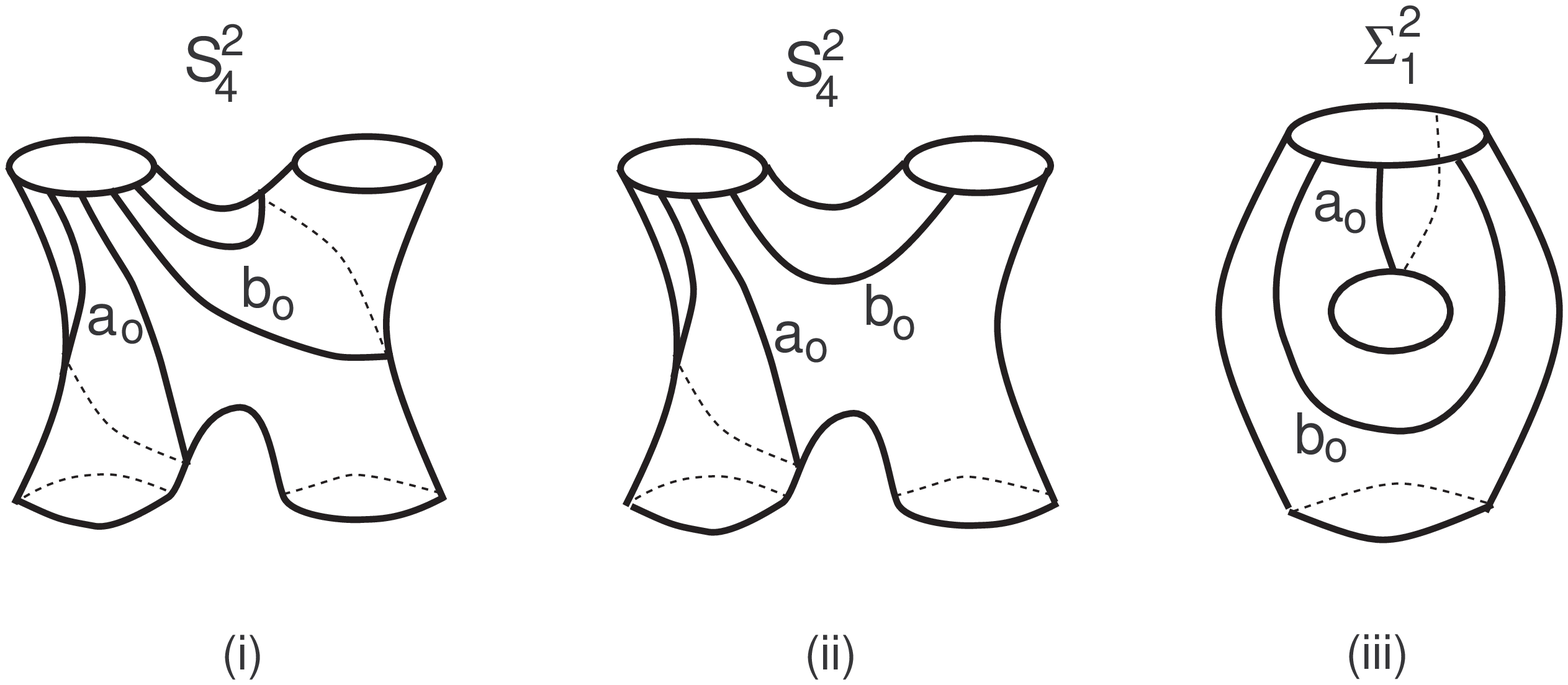}
\caption{Disjoint arcs and neighborhoods}
%Figure 9
\end{center}  \vspace{0.1cm}
\end{figure}

\begin{lemma}
\label{B} Let $a$ and $b$ be two disjoint arcs which are unlinked,
connecting one boundary component $\partial_k$ of $S_c$ to itself
for $k=1,2$. Let $N$ be a regular neighborhood of $a \cup b \cup
\partial_k$ on $S_c$. Then, $(N,a, b) \cong (S_4^2, a_o, b_o)$
where $a_o, b_o$ are the arcs drawn on a standard sphere with four
holes, $S_4^2$, as shown in Figure 9, (i).
\end{lemma}

\begin{lemma}
\label{C} Let $a$ and $b$ be two disjoint arcs on $S_c$ such that
$a$ connects one boundary component $\partial_k$ of $S_c$ to
itself for some $k=1,2$ and $b$ connects the boundary components
$\partial_1$ and $\partial_2$ of $S_c$. Let $N$ be a regular
neighborhood of $a \cup b \cup \partial_1 \cup \partial_2$. Then,
$(N, a, b) \cong (S_4^2, a_o, b_o)$ where $a_o, b_o$ are the arcs
drawn on a standard sphere with four holes, $S_4^2$, as shown in
Figure 9, (ii).
\end{lemma}

\begin{lemma}
\label{D} Let $a$ and $b$ be two disjoint, linked arcs connecting
one boundary component $\partial_k$ of $S_c$ to itself  for
$k=1,2$. Let $N$ be a regular neighborhood of $a \cup b \cup
\partial_k$. Then, $(N, a, b) \cong (\Sigma_1^2, a_o, b_o)$ where
$\Sigma_1^2$ is a standard surface of genus one with two boundary
components, and $a_o, b_o$ are as shown in Figure 9,
(iii).\end{lemma} \vspace{0.1cm}

Let $M$ be a sphere with $k$ holes and $k \geq 5$. A circle $a$ on
$M$ is called an {\it n-circle} if $a$ bounds a disk with $n$
holes on $M$ where $n \geq 2$. If $a$ is a 2-circle on $M$, then
there exists up to isotopy a unique nontrivial embedded arc $a'$
on the two-holed disk component of $M_a$ joining the two holes in
this disc. If $a$ and $b$ are two 2-circles on $M$ such that the
corresponding arcs $a'$, $b'$ can be chosen to meet exactly at one
common end point, and $\alpha = [a], \beta = [b]$, then $\{
\alpha, \beta \}$ is called a {\it simple pair}. A {\it pentagon}
in $\mathcal{C(M)}$ is an ordered 5-tuple $(\alpha_1, \alpha_2,
\alpha_3, \alpha_4, \alpha_5)$, defined up to cyclic permutations,
of vertices of $\mathcal{C(M)}$ such that $i(\alpha_j,
\alpha_{j+1}) = 0$ for $j=1,2,...,5$ and $i(\alpha_j, \alpha_k)
\neq 0$ otherwise, where $\alpha_6 = \alpha_1$. A vertex in
$\mathcal{C(M)}$ is called an {\it n-vertex} if it has a
representative which is an n-circle on $M$. Let $M'$ be the
interior of $M$. There is a natural isomorphism $\chi:
\mathcal{C}(M') \rightarrow \mathcal{C}(M)$ which respects the
above notions and the corresponding notions in \cite{K}. Using
this isomorphism, we can restate a theorem of Korkmaz as follows:

\begin{theorem}
\label{korkmaz} {\bf (Korkmaz) \cite{K}} Let $M$ be a sphere with
$n$ holes and $n \geq 5$. Let $\alpha, \beta$ be two 2-vertices of
$\mathcal{C(M)}$. Then $\{ \alpha, \beta \}$ is a simple pair iff
there exist vertices $\gamma_1, \gamma_2, ..., \gamma_{n-2}$ of
$\mathcal{C}(M)$ satisfying the following conditions: \\
\indent (i) $(\gamma_1, \gamma_2, \alpha, \gamma_3, \beta)$ is a
pentagon
in $\mathcal{C}(M)$, \\
\indent (ii) $\gamma_1$ and $\gamma_{n-2}$ are 2-vertices,
$\gamma_2$ is a 3-vertex and $\gamma_k$ and $\gamma_{n-k}$ are
k-vertices for $3
\leq k \leq \frac{n}{2}$,\\
\indent (iii) $\{\alpha, \gamma_3, \gamma_4, \gamma_5, ...,
\gamma_{n-2}\}$, $\{\alpha, \gamma_2, \gamma_4, \gamma_5, ...,
\gamma_{n-2}\}$,  $\{\beta, \gamma_3, \gamma_4, \gamma_5, ...,
\gamma_{n-2}\}$, and $\{\gamma_1, \gamma_2, \gamma_4, \gamma_5,
..., \gamma_{n-2}\}$ are codimension-zero simplices.\end{theorem}

By using the following lemmas, we will see some more properties of
$\lambda$.

\begin{lemma}
\label{horver} Let $c, x, y, z, h, v$ be essential circles on $S$.
Suppose that there exists a subsurface $N$ of $S$ and a
homeomorphism $\varphi : (N, c, x, y, z, h, v)$ $\rightarrow (N_o,
c_o, x_o, y_o, z_o, h_o, v_o)$ where $N_o$ is a standard sphere
with four holes having $c_o, x_o, y_o, z_o$ on its boundary and
$h_o, v_o$ (horizontal, vertical) are two circles which have
geometric intersection 2 and algebraic intersection 0 as indicated
in Figure 10. Then, there exist $c' \in \lambda([c]), x' \in
\lambda([x]), y' \in \lambda([y]), z' \in \lambda([z]), h' \in
\lambda([h]), v' \in \lambda([v]), N' \subset S$ and a
homeomorphism $\chi: (N', c', x', y', z', h', v') \rightarrow
(N_o, c_o, x_o, y_o, z_o, h_o, v_o)$. \end{lemma}

\begin{figure}[htb]
\begin{center}
\epsfxsize=2.8in \epsfbox{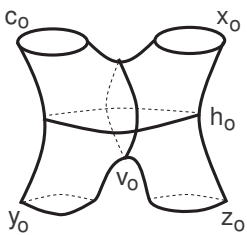} \caption{``Horizontal"
and ``vertical" circles}
\end{center}
\end{figure}

\begin{proof} Let $A$= $\{c, x, y, z, v\}$. Any two elements in $A$
which are isotopic in $S$ bound an annulus on $S$. Let $B$ be a
set consisting of a core from each annulus which is bounded by
elements in $A$, circles in $A$ which are not isotopic to any
other circle in $A$, and $v$. We can extend $B$ to a pants
decomposition $P$ of $S$. Since the genus of $S$ is at least 3,
there are at least four pairs of pants of $P$. Note that $\{v\}$
is a pair of pants decomposition of $N$. Each pair of pants of
this pants decomposition of $N$ is contained in exactly one pair
of pants in $P$. Hence, there is a pair of pants $R$ of $P$ whose
interior is disjoint from $N$ and has at least one of $c, x, y, z$
as one of its boundary components. W.L.O.G. assume that $R$ has
$y$ on its boundary. Let $T$ be a regular neighborhood, in $R$, of
the boundary components of $R$ other than $y$. Let $t, w$ be the
boundary components of $T$ which are in the interior of $R$. Then,
$y, t, w$ bound an embedded pair of pants $Q$ in $R$. Let
$\tilde{N} = N \cup Q$. Then, we can extend $N_o$ to $\tilde{N_o}$
and extend $\varphi$ to a homeomorphism $\tilde{\varphi}:
(\tilde{N}, c, x, y, z, h, v, t, w) \rightarrow (\tilde{N_o}, c_o,
x_o, y_o, z_o, h_o, v_o, t_o, w_o)$, where $\tilde{N_o}$ is as
shown in Figure 11.

Using Lemma \ref{top}, we can choose pairwise disjoint
representatives $c', x', y', z'$, $v', t', w'$ of $\lambda([c]),
\lambda([x]), \lambda([y]), \lambda([z])$, $\lambda([v]),
\lambda([t]), \lambda([w])$ respectively s.t. there exists a
subsurface $\tilde{N'}$ of $S$ and a homeomorphism $\chi :
(\tilde{N'}, c', x', y', z', v', t', w')$ $ \rightarrow
(\tilde{N_o}, c_o, x_o, y_o, z_o, v_o, t_o, w_o)$. Clearly, we
have, $i([h]), [c]) = i([h], [x]) = 0$, $i([h], [y]) = i([h], [z])
= 0$. Since $c, x, y, z$ are all essential circles on $S$, we have
$i([h], [v]) \neq 0$. Then, since $\lambda$ is superinjective, we
have, $i(\lambda([h]), \lambda([c])) =0$, $i(\lambda([h]),
\lambda([x]))= i(\lambda([h]), \lambda([y])) = i(\lambda([h]),
\lambda([z])) = 0$ and $i(\lambda([h]), \lambda([v])) \neq 0$.
Then, a representative $h'$ of $\lambda([h])$ can be chosen such
that $h'$ is transverse to $v'$, $h'$ doesn't intersect any of
$c', x', y', z'$, and $i(\lambda([h]), \lambda([v]) = |h' \cap
v'|$. Since $i(\lambda([h]), \lambda([v])) \neq 0$, $h'$
intersects $v'$. Hence, $h'$ is in the sphere with four holes
bounded by $c', x', y', z'$ in $\tilde{N'}$.

\begin{figure}
\begin{center}
\epsfxsize=2.8in \epsfbox{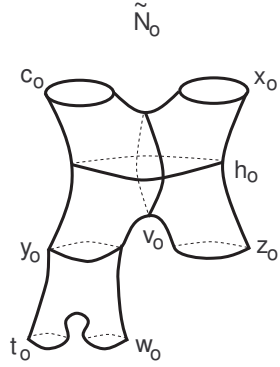} \caption{Sphere with
five holes}
\end{center}
\end{figure}

$\tilde{N}$ and $\tilde{N'}$ are spheres with five holes in $S$.
Since $c, x, y, z$ are essential circles in $S$, the essential
circles on $\tilde{N}$ are essential in $S$. Similarly, since $c',
x', y', z'$ are essential circles in $S$ , the essential circles
on $\tilde{N'}$ are essential in $S$. Furthermore, we can identify
$\mathcal{C}(\tilde{N})$ and $\mathcal{C}(\tilde{N'})$ with two
subcomplexes of $\mathcal{C}(S)$ in such a way that the isotopy
class of an essential circle in $\tilde{N}$ or in $\tilde{N'}$ is
identified with the isotopy class of that circle in $S$. Now,
suppose that $\alpha$ is a vertex in $\mathcal{C}(\tilde{N})$.
Then, with this identification, $\alpha$ is a vertex in
$\mathcal{C}(S)$ and $\alpha$ has a representative in $\tilde{N}$.
Then, $i(\alpha, [c])=i(\alpha, [x])=i(\alpha, [t])=i(\alpha,
[w])= i(\alpha, [z])=0$. Then there are two possibilities: (i)
$i(\alpha, [v]) = 0$ and $i(\alpha, [y]) = 0$ (in this case
$\alpha = [v]$ or $\alpha = [y]$), (ii) $i(\alpha, [v]) \neq 0$ or
$i(\alpha, [y]) \neq 0$. Since $\lambda$ is injective,
$\lambda(\alpha)$ is not equal to any of $[c'], [x'], [t'], [w'],
[z']$. Since $\lambda$ is superinjective, we have,
$i(\lambda(\alpha), [c'])=i(\lambda(\alpha),
[x'])=i(\lambda(\alpha), [t'])=i(\lambda(\alpha), [w'])=
i(\lambda(\alpha), [z'])=0$. Then, there are two possibilities:
(i) $i(\lambda(\alpha), [v']) = 0$ and $i(\lambda(\alpha), [y']) =
0$ (in this case $\lambda(\alpha) = [v']$ or $\lambda(\alpha) =
[y']$), (ii) $i(\lambda(\alpha), [v']) \neq 0$ or
$i(\lambda(\alpha), [y']) \neq 0$. Then, a representative of
$\lambda(\alpha)$ can be chosen in $\tilde{N'}$. Hence, $\lambda$
maps the vertices of $\mathcal{C}(S)$ that have essential
representatives in $\tilde{N}$ to the vertices of $\mathcal{C}(S)$
that have essential representatives in $\tilde{N'}$, (i.e.
$\lambda$ maps $\mathcal{C}(\tilde{N}) \subseteq \mathcal{C}(S)$
to $\mathcal{C}(\tilde{N'}) \subseteq \mathcal{C}(S)$). Similarly,
$\lambda$ maps $\mathcal{C}(N) \subseteq \mathcal{C}(S)$ to
$\mathcal{C}(N') \subseteq \mathcal{C}(S)$.

It is easy to see that $\{[h], [v]\}$ is a simple pair in
$\tilde{N}$. Then, by Theorem \ref{korkmaz}, there exist vertices
$\gamma_1, \gamma_2, \gamma_3$ of $\mathcal{C}(\tilde{N})$ such
that $(\gamma_1, \gamma_2, [h], \gamma_3, [v])$ is a pentagon in
$\mathcal{C}(\tilde{N})$, $\gamma_1$ and $\gamma_3$ are
2-vertices, $\gamma_2$ is a 3-vertex, and $\{[h], \gamma_3\}$,
$\{[h], \gamma_2\}$, $\{[v], \gamma_3\}$ and $\{\gamma_1, \gamma_2
\}$ are codimension-zero simplices of $\mathcal{C}(\tilde{N})$.

Since $\lambda$ is superinjective and $c, x, y, z$ are essential
circles, we can see that $(\lambda(\gamma_1), \lambda(\gamma_2)$,
$\lambda([h]), \lambda(\gamma_3)$, $\lambda([v]))$ is a pentagon
in $\mathcal{C}(\tilde{N'})$. By Lemma \ref{top},
$\lambda(\gamma_1)$ and $\lambda(\gamma_3)$ are 2-vertices, and
$\lambda(\gamma_2)$ is a 3-vertex in $\mathcal{C}(\tilde{N'})$.
Since $\lambda$ is an injective simplicial map $\{\lambda([h]),
\lambda(\gamma_3)\}$, $\{\lambda([h]), \lambda(\gamma_2)\}$,
$\{\lambda([v])$, $\lambda(\gamma_3)\}$ and $\{\lambda(\gamma_1),
\lambda(\gamma_2)\}$ are codimension-zero simplices of
$\mathcal{C}(\tilde{N'})$. Then, by Theorem \ref{korkmaz},
$\{\lambda([h]), \lambda([v])\}$ is a simple pair in $\tilde{N'}$.
Since $\lambda([h])$ has a representative, $h'$, in $N'$, such
that $i(\lambda([h]), \lambda([v]) = |h' \cap v'|$ and
$\{\lambda([h]), \lambda([v])\}$ is a simple pair in $\tilde{N'}$,
there exists a homeomorphism $\chi : (N', c', x', y', z', h', v')
\rightarrow (N_o, c_o, x_o, y_o, z_o, h_o, v_o)$.
\end{proof}

\begin{lemma}
\label{3} Let $c, x, y, z, m, n$ be essential circles on $S$.
Suppose that there exists a subsurface $N$ of $S$ and a
homeomorphism $\varphi: (N, c, x, y, z, m, n) \rightarrow _\varphi
(N_o, c_o, x_o, y_o, z_o, m_o, n_o)$ where $N_o$ is a standard
torus with two boundary components, $c_o, x_o$, and $y_o, z_o,
m_o, n_o$ are circles as shown in Figure 12. Then, there exist $c'
\in \lambda([c]), x' \in \lambda([x]), y' \in \lambda([y]), z' \in
\lambda([z]), m' \in \lambda([m]), n' \in \lambda([n]), N'
\subseteq S$ and a homeomorphism $\chi$ such that $(N', c', x',
y', z', m', n') \rightarrow _\chi (N_o, c_o,$ $x_o, y_o, z_o, m_o,
n_o)$.

\begin{figure}[htb]
\begin{center}
\hspace{1.3cm} \epsfxsize=3.7in \epsfbox{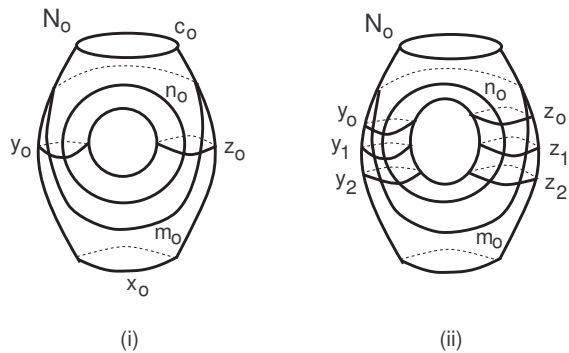}
\caption{Circles on torus with two boundary components}
%Figure 12
\end{center}
\end{figure}
\end{lemma}

\begin{proof} Since the genus of $S$ is at least 3, $c$ cannot
be isotopic to $x$ in $S$. So, we can complete $\{c, x, y, z\}$ to
a pair of pants decomposition, $P$, of $S$. Since $\{y, z\}$ gives
a pair of pants decomposition on $N$, by Lemma \ref{top}, there
exists a subsurface $N' \subseteq S$ which is homeomorphic to
$N_o$ and there are pairwise disjoint representatives $c', x', y',
z'$ of $\lambda([c]), \lambda([x]), \lambda([y]), \lambda([z])$
respectively and a homeomorphism $\phi$ such that $(N', c', x',
y', z') \rightarrow _\phi (N_o, c_o, x_o, y_o, z_o)$. Then by
Lemma \ref{intone}, we have the following:

\noindent $i([m], [z])=1 \Rightarrow i(\lambda([m]),
\lambda([z]))=1$, $i([n], [z])=1 \Rightarrow i(\lambda([n]),
\lambda([z]))=1$, $i([m], [y])=1 \Rightarrow i(\lambda([m]),
\lambda([y]))=1$, $i([n], [y])=1 \Rightarrow i(\lambda([n]),
\lambda([y]))=1$, $i([m], [c])=0 \Rightarrow i(\lambda([m]),
\lambda([c]))=0$, $i([n], [c])=0 \Rightarrow i(\lambda([n]),
\lambda([c]))=0$, $i([m], [x])=0 \Rightarrow i(\lambda([m]),
\lambda([x]))=0$, $i([n], [x])=0 \Rightarrow i(\lambda([n]),
\lambda([x]))=0$, $i([m], [n])=0 \Rightarrow i(\lambda([m]),
\lambda([n]))=0$.

There are representatives $m_1 \in \lambda([m]), n' \in
\lambda([n])$ such that $|m_1 \cap y'| = |m_1 \cap z'| = 1$, $|m_1
\cap c'| = |m_1 \cap x'| = 0$, $|n' \cap y'| = |n' \cap z'| = 1,
|n' \cap c'| = |n' \cap x'| = |m_1 \cap n'| = 0$ with all
intersections transverse.

Since $\phi$ is a homeomorphism, we have $|\phi(m_1) \cap
\phi(n')|= 0, |\phi(n') \cap y_o| = 1$, $|\phi(n') \cap z_o| = 1,
|\phi(n') \cap c_o| = |\phi(n') \cap x_o| = 0 = |\phi(m_1) \cap
c_o| = |\phi(m_1) \cap x_o| = 0$, $|\phi(m_1) \cap y_o| =
|\phi(m_1) \cap z_o| =1$.

Let's choose parallel copies $y_1, y_2$ of $y_o$ and $z_1, z_2$ of
$z_o$ as shown in Figure 12 so that each of them has transverse
intersection one with $\phi(m')$ and $\phi(n')$. Let $P_1, P_2$ be
the pair of pants with boundary components $c_o, y_o, z_o$, and
$x_o, y_2, z_2$ respectively. Let $Q_1, Q_2, R_1, R_2$ be the
annulus with boundary components $\{y_o, y_1\}, \{y_1, y_2\}$,
$\{z_o, z_1\}$, $\{z_1, z_2\}$ respectively. By the classification
of isotopy classes of families of properly embedded disjoint arcs
in pairs of pants, $\phi(m') \cap P_1$, $\phi(m') \cap P_2$,
$\phi(n') \cap P_1$ and $\phi(n') \cap P_2$ can be isotoped to the
arcs $m_o \cap P_1, m_o \cap P_2, n_o \cap P_1, n_o \cap P_2$
respectively. Let $\kappa : P_1 \times I \rightarrow P_1$, $\tau :
P_2 \times I \rightarrow P_2$ be such isotopies. By a tapering
argument, we can extend $\kappa $ and $\tau $ and get
$\tilde{\kappa }: (P_1 \cup Q_1 \cup R_1) \times I \rightarrow
(P_1 \cup Q_1 \cup R_1)$ and $\tilde{\tau }: (P_2 \cup Q_2 \cup
R_2) \times I \rightarrow (P_2 \cup Q_2 \cup R_2)$ so that
$\tilde{\kappa }_t$ is $id$ on $y_1 \cup z_1$ and $\tilde{\tau
}_t$ is $id$ on $y_1 \cup z_1$ for all $t \in I$. Then, by gluing
these extensions we get an isotopy $\vartheta$ on $N_o \times I$.

By the classification of isotopy classes of arcs (relative to the boundary) on an annulus, $\vartheta_1(\phi(n'))
\cap (R_1 \cup R_2)$ can be isotoped to $t_{z_o} ^k(n_o) \cap (R_1 \cup R_2)$ for some $k \in \mathbb{Z}$.
Let's call this isotopy $\mu$. Let $\tilde{\mu }$ denote the extension by id to $N_o$. Similarly, $\vartheta_1
(\phi(n')) \cap (Q_1 \cup Q_2)$ can be isotoped to $t_{y_o} ^l (n_o) \cap (Q_1 \cup Q_2)$ for some $l \in
\mathbb{Z}$. Let's call this isotopy $\nu $. Let $\tilde{\nu }$ denote the extension by id to $N_o$. Then,
``gluing'' the two isotopies $\tilde{\mu}$  and $\tilde{\nu}$, we get a new isotopy, $\epsilon$, on $N_o$.
Then we have, $t_{y_o} ^{-l} (t_{z_o} ^{-k}(\epsilon _1(\vartheta_1 (\phi(n'))))) = n_o$.
%So, $\phi^{-1} (\epsilon_1 ^{-1} (t_{z_o} ^k (n_o))) = n'$.
Clearly, $t_{y_o} ^{-l} \circ t_{z_o} ^{-k} \circ \epsilon_1 \circ
\vartheta _1$ fixes $c_o, x_o, y_o, z_o$. So, we get $t_{y_o}
^{-l} \circ t_{z_o}^{-k} \circ \epsilon_1 \circ \vartheta_1 \circ
\phi : (N', c', x', y', z', n') \rightarrow (N_o, c_o, x_o, y_o,
z_o, n_o)$. Let $\chi = t_{y_o} ^{-l} \circ t_{z_o}^{-k} \circ
\epsilon_1 \circ \vartheta_1 \circ \phi$. Then, we also get,
$\chi(m_1)$ is isotopic to $m_o$ because of the intersection
information. Let $m'= \chi ^{-1}(m_o)$. Then we get, $\chi : (N',
c', x', y', z', m', n') \rightarrow (N_o, c_o, x_o, y_o, z_o, m_o,
n_o)$.\end{proof}

Let $i$ be an essential properly embedded arc on $S_c$. Let $A$ be
a boundary component of $S_c$ which has one end point of $i$ and
$B$ be the boundary component of $S_c$ which has the other end
point of $i$. Let $N$ be a regular neighborhood of $i \cup A \cup
B$ in $S_c$. By Euler characteristic arguments, $N$ is a pair of
pants. The boundary components of $N$ are called {\it encoding
circles of $i$ on $S_c$}. These circles correspond to nontrivial
circles on $S$, which are called {\it encoding circles of $i$ on
$S$}. The set of isotopy classes of encoding circles on $S$ is
called the {\it encoding simplex}, $\Delta_i$, of $i$ (and of
$[i]$). Note that $[c]$ is a vertex in the encoding simplex of
$i$.

An essential properly embedded arc $i$ on $S_c$ is called {\it
type 1.1} if it joins one boundary component $\partial_k$ of $S_c$
to itself for $k=1,2$ and if $\partial_k \cup i$ has a regular
neighborhood $N$ in $S_c$ which has only one circle on its
boundary which is inessential w.r.t. $S_c$. If $N$ has two circles
on its boundary which are inessential w.r.t. $S_c$, then $i$ is
called {\it type 1.2}. We call $i$ to be {\it type 2}, if it joins
the two boundary components $\partial_1$ and $\partial_2$ of $S_c$
to each other. An element $[i] \in \mathcal{V}(S_c)$ is called
{\it type 1.1 (1.2, 2)} if it has a type 1.1 (1.2, 2)
representative.

Let $\partial^1, \partial^2$ be the boundary components of $S_d$.
We prove the following lemmas in order to show that $\lambda$
induces a map $\lambda_* : \mathcal{V}(S_c) \rightarrow
\mathcal{V}(S_d)$ with certain properties.

\begin{lemma}
\label{boundary} Let $\partial_k \in \partial S_c$ for some $k \in
\{1,2\}$. Then, there exists a unique $\partial^l \in \partial
S_d$ for some $l \in \{1,2 \}$ such that if $i$ is a properly
embedded essential arc on $S_c$ connecting $\partial_k$ to itself,
then there exists a properly embedded arc $j$ on $S_d$ connecting
$\partial^l$ to itself such that $\lambda(\Delta_i) = \Delta_j$.
\end{lemma}

\begin{proof} Assume that each of $\partial^1$ and
$\partial^2$ satisfies the hypothesis. Let $i$ be a properly
embedded, essential, type 1.1 arc connecting $\partial_k$ to
itself. Then, there exist properly embedded arcs, $j_1$,
connecting $\partial^1$ to itself, and $j_2$, connecting
$\partial^2$ to itself, such that $\lambda(\Delta_i) =
\Delta_{j_1}$ and $\lambda(\Delta_i) = \Delta_{j_2}$. Then, we
have $\Delta_{j_1} = \Delta_{j_2}$. Note that a properly embedded
essential arc $i$ is type 1.1 iff $\Delta_i$ has exactly 3
elements. Otherwise $\Delta_i$ has 2 elements. Since $i$ is type
1.1 and $\lambda(\Delta_i) = \Delta_{j_1}$ and $\lambda(\Delta_i)
= \Delta_{j_2}$ and $\lambda$ is injective, $j_1$ and $j_2$ are
type 1.1. We can choose a pairwise disjoint representative set
$\{a, b, d\}$ of $\Delta_{j_1}$ on $S$. Since $\Delta_{j_1} =
\Delta_{j_2}$, $\{a, b, d\}$ is a pairwise disjoint representative
set for $\Delta_{j_2}$ on $S$. Then, the curves $\tilde{a},
\tilde{b}$ on $S_d$ which correspond to $a, b$ on $S$ and
$\partial^1$ bound a pair of pants, $P$, on $S_d$ containing an
arc, $j_1'$, isotopic to $j_1$. Similarly, $\tilde{a}, \tilde{b},
\partial^2$ bound a pair of pants, $Q$, on $S_d$ containing an
arc, $j_2'$, isotopic to $j_2$. Let's cut $S_d$ along $\tilde{a}$
and $\tilde{b}$. Then, $P$ is the connected component of $S_{d
\cup a \cup b}$ containing $\partial^1$ and $Q$ is the connected
component of $S_{d \cup a \cup b}$ containing $\partial^2$. $P
\neq Q$ since $\partial^2$ is not in $P$ and $\partial^2$ is in
$Q$. Then $P$ and $Q$ are distinct connected components meeting
along $\tilde{a}$ and $\tilde{b}$. Hence, $S_d$ is $P \cup Q$, a
torus with two holes. This implies that $S$ is a genus 2 surface
which gives a contradiction since the genus of $S$ is at least 3.
So, only one boundary component of $S_d$ can satisfy the
hypothesis.

Since $i$ is type 1.1, $\Delta_i$ contains $[c]$ and two other
isotopy classes of nontrivial circles which are not isotopic to
$c$ in $S$. Let $P'$ be a pairwise disjoint representative set of
$\lambda([\Delta_i])$, containing $d$. By the proof of Lemma
\ref{top}, $P'$ bounds a pair of pants on $S$. Since the genus of
$S$ is at least 3, $P'$ bounds a unique pair of pants on $S$,
which corresponds to a unique pair of pants, $Q$, in $S_d$ which
has only one inessential boundary component. Let $\partial^{l(i)}$
be this inessential boundary component. Let $j$ be an essential
properly embedded arc connecting $\partial^{l(i)}$ to itself in
$Q$. Then, we have $\lambda(\Delta_i) = \Delta_j$.

Now, to see that $\partial^{l(i)}$ is independent of the type 1.1
arc $i$ connecting $\partial_k$ to itself, we prove the following
claim:

Claim I: If we start with two type 1.1 arcs $i$ and $j$ starting
and ending on $\partial_k$, then $\partial^{l(i)} =
\partial^{l(j)}$.

Proof of Claim I: Let $[i], [j]$ be type 1.1 and $i, j$ connect
$\partial_k$ to itself. W.L.O.G. we can assume that $i$ and $j$
have minimal intersection. First, we prove that there is a
sequence $j = r_0 \rightarrow r_1 \rightarrow ... \rightarrow
r_{n+1}=i$ of essential properly embedded arcs joining
$\partial_k$ to itself so that each consecutive pair is disjoint,
i.e. the isotopy classes of these arcs define a path in
$\mathcal{B} (S_c)$, between $i$ and $j$.

If $|i \cap j|=0$, then take $r_0=j$, $r_1=i$. We are done. Assume
that $|i \cap j|=m>0$. We orient $i$ and $j$ arbitrarily. Then, we
define two arcs in the following way: Start on the boundary
component $\partial_k$, on one side of the beginning point of $j$
and continue along $j$ without intersecting $j$, till the last
intersection point along $i$. Then we would like to follow $i$,
without intersecting $j$, until we reach $\partial_k$. So, if we
are on the correct side of $j$ we do this; if not, we change our
starting side from the beginning and follow the construction. This
gives us an arc, say $j_1$. We define $j_2$, another arc, by
changing the orientation of $j$ and following the same
construction. It is easy to see that $j_1, j_2$ are disjoint
properly embedded arcs connecting $\partial_k$ to itself as $i$
and $j$ do. One can see that $j_1, j_2$ are essential arcs since
$i, j$ intersect minimally. In Figure 13, we show the beginning
and the end points of $i$, the essential intersections of $i, j$,
and $j_1, j_2$ near the end point of $i$ on $\partial_k$.

\begin{figure}[htb]
\begin{center}
\epsfxsize=2.65in \epsfbox{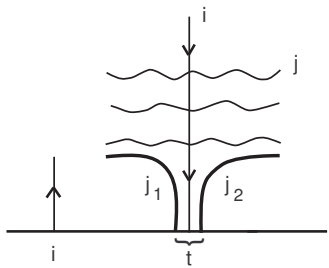} \caption{Splitting the
arc $j$ along the end of $i$}
%Figure 13
\end{center}
\end{figure}

$|i \cap j_1| < m$, $|i \cap j_2| < m$ since we eliminated at
least one intersection with $i$. We also have $|j_1 \cap j| = |j_2
\cap j| =0$ since we never intersected $j$ in the construction.
Notice that $j_1$ and $j_2$ are not oriented, and $i$ is oriented.

Claim II: $j_1$ is not isotopic to $j_2$.

Proof of Claim II: Suppose that $j_1$ and $j_2$ are isotopic.
Then, they are not linked, and there is a band $B$ such that
$\partial B \subseteq j_1 \cup j_2 \cup \partial_k$ and $\partial
B \setminus (j_1 \cup j_2)$ is a disjoint union of two arcs
$\varepsilon_1, \varepsilon_2$ on $\partial_k$, where each of
$\varepsilon_1$ and $\varepsilon_2$ starts at an end point of
$j_1$ and ends at an end point of $j_2$. By the construction,
there is an end point of each of $j_1$ and $j_2$ near the end
point of $i$ on $\partial_k$, on the arc $t$ shown in the figure,
and there is not any point of $j$ on $t$. Then, either $t =
\varepsilon_1$ or $t = \varepsilon_2$. W.L.O.G. assume that $t=
\varepsilon_1$. Then, $\varepsilon_1 \subseteq \partial B$ has the
end point of $i$. Since $\varepsilon_1 \subseteq \partial_k$, it
has only one side on $S_c$, and hence on $B$. Then, since
$\varepsilon_1 \subseteq
\partial B$ has the end point of $i$, by the construction we can
see that the last intersection point along $i$ of $i$ and $j$ has
to lie in the band $B$. Then, since $j$ does not intersect any of
$j_1$ and $j_2$, and $B$ has a point of $j$ in the interior, $j$
has to live in $B$. Then, since $\varepsilon_1$ doesn't contain
any point of $j$, and $j$ does not intersect any of $j_1$ and
$j_2$, the end points of $j$ has to lie on $\varepsilon_2$, which
implies that $j$ is an inessential arc on $S_c$. This gives a
contradiction. Hence, $j_1$ and $j_2$ cannot be isotopic.

Claim III: Either $j_1$ or $j_2$ is of type 1.1.

Proof of Claim III: If $j_1$ and $j_2$ are linked, then a regular
neighborhood of $j_1 \cup j_2 \cup \partial_k$ in $S_c$ is a genus
one surface with two boundary components by Lemma \ref{D}. In this
case, both arcs have to be type 1.1 since encoding circles on
$S_c$ for $j_1, j_2$ can be chosen in this surface and this
surface has only one boundary component which is a boundary
component of $S_c$. So, both arcs have to be type 1.1.

If $j_1$ and $j_2$ are unlinked, then a regular neighborhood, $N$,
of $j_1 \cup j_2 \cup \partial_i$ in $S_c$ is a sphere with four
holes by Lemma \ref{B}. Encoding circles of $j_1, j_2$ on $S_c$
can be chosen in $N$. Assume that $j_1$ is type 1.2 and has two
encoding circles which are the two boundary components of the
surface $S_c$. Then, $j_2$ would have to be type 1.1 since $N$ can
have at most two boundary components which are the boundary
components of $S_c$ and regular neighborhoods of $j_1 \cup
\partial_k$ and $j_2 \cup
\partial_k$ have only one common boundary component which is
$\partial_k$. This proves Claim III.

Let $r_1 \in \{j_1, j_2\}$ and $r_1$ be type 1.1. By the
construction we get, $|i \cap r_1| < m$, $|j \cap r_1|=0$. Now,
using $i$ and $r_1$ in place of $i$ and $j$ we can define a new
type 1.1 arc $r_2$ with the properties, $|i \cap r_2| < |i \cap
r_1|, |r_1 \cap r_2|=0$. By an inductive argument, we get a
sequence of arcs such that every consecutive pair is disjoint,
$i=r_{n+1} \rightarrow r_n \rightarrow r_{n-1} \rightarrow ...
\rightarrow r_1 \rightarrow r_0=j$. This gives us a path of type
1.1 arcs in $\mathcal{B}(S_c)$ between $i$ and $j$. By using Lemma
\ref{B} and Lemma \ref{D}, we can see a regular neighborhood of
the union of each consecutive pair in the sequence and the
boundary components of $S_c$, and encoding circles of these
consecutive arcs. Then, by using the results of Lemma \ref{horver}
and \ref{3}, we can see that each pair of disjoint type 1.1 arcs
give us the same boundary component. Hence, by using the sequence
given above, we conclude that $i$ and $j$ give us the same
boundary component. This proves Claim I.

Let $i_o$ be a properly embedded, type 1.1 arc on $S_c$ connecting
$\partial_k$ to itself. Let $\partial^l= \partial^{l(i_o)}$. If
$i$ is a properly embedded, type 1.1 arc on $S_c$ connecting
$\partial_k$ to itself, then by the arguments given above we have
$\partial^l = \partial^ {l(i)}$, and there exists a properly
embedded arc $j$ on $S_d$ connecting $\partial^l$ to itself such
that $\lambda (\Delta_i) = \Delta_j$. Suppose that $i$ is a
properly embedded, type 1.2 arc on $S_c$ connecting $\partial_k$
to itself. Let $\Delta_i$ be the encoding simplex of $[i]$. Since
$i$ is type 1.2, $\Delta_i$ contains $[c]$ and only one other
isotopy class of a nontrivial circle which is not isotopic to $c$
in $S$. A pairwise disjoint representative set, $P$, of $\Delta_i$
corresponds to a nonembedded pair of pants on $S$. Let $P'$ be a
pairwise disjoint representative set of $\lambda([\Delta_i])$,
containing $d$. By extending $P$ to a pair of pants decomposition
of $S$, and applying Lemma \ref{top}, we can see that $P'$
corresponds to a unique nonembedded pair of pants of $S$ which
corresponds to a unique pair of pants, $Q$, in $S_d$ which has two
inessential boundary components containing $\partial^l$. Let $j$
be a nontrivial properly embedded arc connecting $\partial^l$ to
itself in $Q$. Then, we have $\lambda(\Delta_i) = \Delta_j$.
Hence, $\partial^l$ is the boundary component that we
want.\end{proof}

We define a map $\sigma: \{ \partial_1, \partial_2 \} \rightarrow
\{ \partial^1, \partial^2 \}$ using the correspondence which is
given by Lemma \ref{boundary}.

\begin{lemma} $\sigma$ is a bijection.
\end{lemma}

\begin{proof}
Let $x, y, z, t, h, c$ be essential circles on $S$ such that $x,
y, z, t$ bound a 4-holed sphere and $h, c$ be two circles which
intersect geometrically twice and algebraically zero times in this
subsurface as shown in Figure 14, (i). By using Lemma
\ref{horver}, we can see that there exist pairwise disjoint
representatives $x', y', z', t', h', d$ of $\lambda([x]),
\lambda([y])$, $\lambda([z]), \lambda([t]), \lambda([h]),
\lambda([c])$ respectively such that $x', y', z', t'$ bound a
sphere with four holes, and $h', d$ intersect geometrically twice
and algebraically zero times in this subsurface as shown in Figure
14, (ii).

\begin{figure}[htb]
\begin{center}
\hspace{1.3cm} \epsfxsize=3.6in \epsfbox{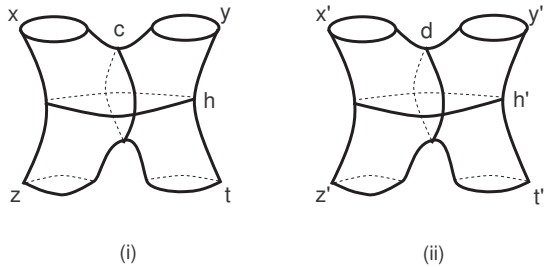}
\caption{Sphere with four holes}
%Figure 13
\end{center}
\end{figure}

The curve $h$ corresponds to two arcs, say $i_1, i_2$, on $S_c$
which start and end on different boundary components,
$\partial_k$, $\partial_l$, of $S_c$ respectively. W.L.O.G. assume
that $\{x, z, c\}$ and $\{y, t, c\}$ are encoding circles of $i_1$
and $i_2$ on $S$. $\{x, z, c\}$ and $\{y, t, c\}$ correspond to
$\{x', z', d\}$ and $\{y', t', d\}$ which bound pairs of pants $L$
and $R$ on $S$ respectively. Assume that $L$ corresponds to a pair
of pants, $L'$, on $S_d$ which has a boundary component
$\partial^m$ of $S_d$ and $R$ corresponds to a pair of pants,
$R'$, on $S_d$ which has a boundary component $\partial^n$ of
$S_d$. $h'$ corresponds to two arcs, say $j_1$ and $j_2$ on $L'$
and $R'$ which start and end on different boundary components,
$\partial^m$ and $\partial^n$ of $S_d$ respectively. By the
definition of $\sigma$ we have, $\sigma(\partial_k) =
\partial^m$ and $\sigma(\partial_l) = \partial^n$ for $k, l, m, n
\in \{1,2\}$, $k \neq l, m \neq n$. So, $\sigma$ is onto. Hence,
$\sigma$ is a bijection.\end{proof}

\begin{lemma}
\label{vertex} Let $[i] \in \mathcal{V}(S_c)$. If $i$ connects
$\partial_k$ to $\partial_l$ on $S_c$ where $k, l \in \{1,2\}$,
then there exists a unique $[j] \in  \mathcal{V}(S_d)$ such that
$j$ connects $\sigma(\partial_k)$ to $\sigma(\partial_l)$ and
$\lambda(\Delta_i) = \Delta_j$.
\end{lemma}

\begin{proof} Let $[i] \in \mathcal{V}(S_c)$ and let $i$ connect
$\partial_k$ to $\partial_l$ on $S_c$ where $k, l \in \{1,2\}$.
Let $\Delta_i$ be the encoding simplex of $i$. Then, a pairwise
disjoint representative set of $\lambda([\Delta_i])$ containing
$d$ corresponds to a unique pair of pants in $S_d$ which has
boundary components $\sigma(\partial_k)$ and $\sigma(\partial_l)$.
By the classification of properly embedded arcs in pair of pants,
there exists a unique isotopy class of nontrivial properly
embedded arcs which start at $\sigma(\partial_k)$ and end at
$\sigma(\partial_l)$, in this pair of pants. Let $j$ be such an
arc. We have $\lambda(\Delta_i) = \Delta_j$.

Let $e$ be an essential properly embedded arc in $S_d$ such that
$e$ connects $\sigma(\partial_k)$ to $\sigma(\partial_l)$ and
$\lambda(\Delta_i) = \Delta_e$. Then, we have $\Delta_e = \Delta_j
= \lambda(\Delta_i)$. Let $a, b, d$ be a pairwise disjoint
representative set of $\lambda(\Delta_i)$ on $S$. Then there are
properly embedded arcs $j_1, e_1$ isotopic to $j, e$ respectively
which are in pair of pants bounded by $a, b, d$. Since the genus
of $S$ is at least 3, there is at most one pair of pants which has
$a, b, d$ on its boundary. So, $j_1, e_1$ are in the same pair of
pants. Since they both connect the same boundary components in
this pair of pants, they are isotopic. So, $[j] = [e]$. Hence,
$[j]$ is the unique isotopy class in $S_d$ such that $j$ connects
$\sigma(\partial_k)$ to $\sigma(\partial_l)$ and
$\lambda(\Delta_i) = \Delta_j$. \end{proof}

$\lambda$ induces a unique map $\lambda_* : \mathcal{V}(S_c)
\rightarrow \mathcal{V}(S_d)$ such that if $[i] \in
\mathcal{V}(S_c)$ then $\lambda_*([i])$ is the unique isotopy
class corresponding to $[i]$ where the correspondence is given by
Lemma \ref{vertex}. Using the results of the following lemmas, we
will prove that $\lambda_*$ extends to an injective simplicial map
on the whole complex, $\mathcal{B}(S_c)$.

\begin{lemma} $\lambda_* :\mathcal{V}(S_c) \rightarrow \mathcal{V}(S_d)$
extends to a simplicial map
$\lambda_* : \mathcal{B}(S_c) \rightarrow \mathcal{B}(S_d)$.
\end{lemma}

\begin{proof} It is enough to prove that if two distinct isotopy
classes of essential properly embedded arcs in $S_c$ have disjoint
representatives, then their images under $\lambda_*$ have disjoint
representatives. Let $a, b$ be two disjoint representatives of two
distinct classes in $\mathcal{V}(S_c)$. Let $\partial_1,
\partial_2$ be the two boundary components of $S_c$. We consider
the following cases:

\begin{figure}[htb]
\begin{center}
\epsfxsize=2.53in \epsfbox{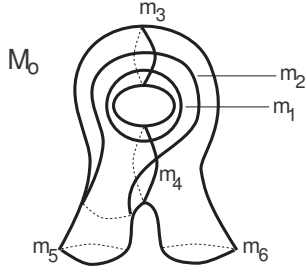} \caption{Circles
associated to type 2 arcs}
\end{center}
\end{figure}

Case 1: Assume that $a$ and $b$ connect the two boundary
components of $S_c$. W.L.O.G. assume that the end points of $a$
map to the same point under the quotient map $q: S_c \rightarrow
S$, and the end points of $b$ map to the same point under $q$. By
Lemma \ref{A}, we have, $(N, a, b) \cong (S^2_4, a_o, b_o)$ where
$N$ is a regular neighborhood of $a \cup b \cup
\partial_1 \cup \partial_2$. Then, we have $(q(N), q(a), q(b),
q(\partial_1)) \cong (M_o, m_1, m_2, m_3)$ where $M_o, m_1, m_2,
m_3$ are as shown in Figure 15.

It is easy to see that $(M_o, m_1, m_2, m_3) \cong (N_o, n_o, m_o,
z_o)$, where $N_o, n_o, m_o$, $z_o$ are as given in Lemma \ref{3}.
Then, there exists a homeomorphism $\varphi: (q(N), q(a)$, $q(b),
q(\partial_1)) \rightarrow (N_o, n_o, m_o, z_o)$. Let $v =
\varphi^{-1}(m_4), y = \varphi^{-1}(m_5), z = \varphi^{-1}(m_6)$
where $m_4, m_5, m_6$ are as shown in Figure 15. Suppose that
$q(a)$ is isotopic to $q(b)$ on $S$. Since $q(a)$ is disjoint from
$q(b)$, there exists an annulus having $q(a) \cup q(b)$ as its
boundary. Since each of $q(a)$ and $q(b)$ intersects $c$
transversely once, $c$ cuts this annulus into a band $B$ on $S_c$
such that $\partial B \subseteq a \cup b \cup \partial S_c$. But
this is not possible since $a$ is not isotopic to $b$ on $S_c$.
Hence, $q(a)$ is not isotopic to $q(b)$ on $S$. Then, $y$ and $z$
are essential circles on $S$. Since $q(\partial_1)$ is equal to
$c$ on $S$, $q(\partial_1)$ is an essential circle on $S$. Since
each of $q(a)$ and $q(b)$ intersects $c$ transversely once, $q(a)$
and $q(b)$ are essential circles on $S$. Since $v$ intersects
$q(b)$ transversely once, $v$ is also essential on $S$.  Note that
$d \in \lambda([q(\partial_1)])$. Then, by Lemma \ref{3}, we can
find $a' \in \lambda([q(a)])$, $b' \in \lambda([q(b)])$, $N'
\subseteq S$ and a homeomorphism $\chi : (N', a', b', d)
\rightarrow (N, n_o, m_o, z_o)$. If we cut $N'$ along $d$, we get
two disjoint arcs, $\hat{a}, \hat{b}$ corresponding to $a', b'$
respectively.

Since $a$ connects $\partial_1$ and $\partial_2$, the encoding
simplex of $[a]$ is $\{[c], \gamma\}$ where $\gamma$ is the class
of a 1-separating circle on $S$. There exists $x \in \gamma$ such
that $x$ bounds a genus 1 subsurface $Q$ and $q(a)$ and $c$
intersect transversely once on $Q$. Since $i([x], [c])=0$, $i([x],
[q(a)]) =0$, and $\lambda$ is superinjective, we have $i(
\lambda([x]), [d])=0$, $i(\lambda([x]), \lambda([q(a)]))$ $= 0$.
Then we choose a representative $x'$ of $\lambda([x])$ which is
disjoint from $a' \cup d$. By Lemma 3.7, $x'$ is a 1-separating
circle bounding a subsurface $R$ containing $a'\cup d$. Then, it
is easy to see that $\{[x'], [d]\}$ is the encoding simplex of
$[\hat{a}]$. Then, since the encoding simplex of $[\hat{a}]$ is
the image of the encoding simplex of $[a]$ under $\lambda$, and
both $a$ and $[\hat{a}]$ are type 2, by the definition of
$\lambda_*$, $[\hat{a}]$ is the image of $[a]$ under $\lambda_*$.
Similarly, $[\hat{b}]$ is the image of $[b]$ under $\lambda_*$.
This shows that $\lambda_*([a])$ and $\lambda_*([b])$ have
disjoint representatives $\hat{a}$ and $\hat{b}$ respectively.

Case 2: Assume that $a$ connects one boundary component of $S_c$
to itself and $b$ connects the other boundary component of $S_c$
to itself. Then, since $a$ and $b$ are disjoint, they are both
type 1.1 arcs. W.L.O.G. we can assume that $a$ connects
$\partial_1$ to itself and $b$ connects $\partial_2$ itself in
such a way that $q(\partial a) = q(\partial b)$ where $q$ is the
quotient map $q: S_c \rightarrow S$. Then, since $a \cup
\partial_1$ is disjoint from $b \cup \partial_2$, we can find
disjoint regular neighborhoods, $N_1$, of $a \cup \partial_1$ and
$N_2$ of $b \cup \partial_2$ on $S_c$. Then, $q(N_1 \cup N_2)$ is
a 4-holed sphere. Let $x, y, z, t$ be the boundary components of
this subsurface s.t. $x, z, c$ are encoding circles for $a$ on $S$
and $y, t, c$ are encoding circles for $b$ on $S$. Let $h = q(a
\cup b)$. Then, $h$ and $c$ intersect geometrically twice and
algebraically zero times in this subsurface as shown in Figure 14,
(i).

By using Lemma \ref{horver}, we can see that there exist pairwise
disjoint representatives $x', y', z', t', h', d$ of $\lambda([x]),
\lambda([y])$, $\lambda([z]), \lambda([t]), \lambda([h]),
\lambda([c])$ respectively such that $x', y', z', t'$ bound a
4-holed sphere, and $h', d$ intersect geometrically twice and
algebraically zero times in this subsurface as shown in Figure 14,
(ii).

The curve $h'$ corresponds to two essential disjoint arcs when we
cut this four holed sphere along $d$. As in Case 1, by the
definition of $\lambda_*$, these two disjoint arcs are
representatives for $\lambda_*([a])$ and $(\lambda_*([b])$ which
proves the lemma for Case 2.

\begin{figure}[htb]
\begin{center}
\epsfxsize=3in \epsfbox{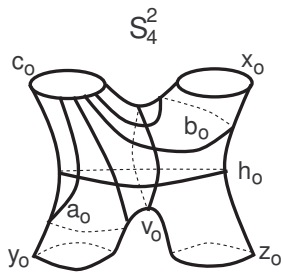} \caption{Arcs and their
encoding circles, I}
%Figure 16
\end{center}
\end{figure}

Case 3: Assume that $a, b$ are unlinked, connecting $\partial_i$
to itself for some $i=1,2.$ By Lemma \ref{B}, there is a
homeomorphism $\phi$ such that $(S^2_4, a_o, b_o) \cong_\phi (N,
a, b)$ where $N$ is a regular neighborhood of $a \cup b \cup
\partial_i$ and $a_o, b_o$ are as shown in Figure 16.

Since $h_o, x_o, c_o$ and $v_o, y_o, c_o$ are the boundary
components of regular neighborhoods of $c_o \cup b_o$ and $c_o
\cup a_o$ on $S_4^2$ respectively, $\phi(h_o)$, $\phi(x_o)$ and
$\phi(v_o), \phi(y_o)$ are encoding circles for $b$ and $a$ on
$S_c$ respectively. Then, $q(\phi(v_o)), q(\phi(h_o))$ are
essential circles on $S$ where $q: S_c \rightarrow S$ is the
quotient map. We have, $(S^2_4, c_o, x_o, y_o, z_o, h_o, v_o)
\cong (N, \phi(c_o), \phi(x_o), \phi(y_o), \phi(z_o), \phi(h_o),
\phi(v_o))$. Since $a$ and $b$ are essential arcs on $S_c$,
$q(\phi(x_o))$ and $q(\phi(y_o))$ are essential circles on $S$.
$q(\phi(c_o))$ is equal to $c$ on $S$. So, $q(\phi(c_o))$ is an
essential circle on $S$. Since $a$ is not isotopic to $b$ on
$S_c$, $q(\phi(z_o))$ is an essential circle on $S$. So, we have
that all of $q(\phi(v_o)), q(\phi(h_o))$, $q(\phi(c_o)),
q(\phi(x_o)), q(\phi(y_o)), q(\phi(z_o))$, are essential circles
on $S$. Note that $d \in \lambda([\phi(c_o)])$. Then, by using
Lemma \ref{horver}, we can find $x' \in \lambda ([\phi (x_o)]), y'
\in \lambda ([\phi (y_o)])$, $z' \in \lambda ([\phi (z_o)]), h'
\in \lambda ([\phi (h_o)]), v' \in \lambda ([\phi(v_o)])$, $N'
\subseteq S$ and a homeomorphism $\chi : (S^2_4, c_o, x_o, y_o,
z_o, h_o, v_o) \rightarrow (N', d, x', y', z', h', v')$. Then,
$h', x'$ are encoding circles for $\lambda_*([b])$ on $S$ and $v',
y'$ are encoding circles for $\lambda_*([a])$ on $S$, and $\chi(
b_o ), \chi(a_o)$ are disjoint representatives for $\lambda_*
([b])$ and $\lambda_* ([a])$ respectively.

Case 4: Assume that $a$ connects one boundary component of $S_c$
to itself and $b$ connects the two boundary components of $S_c$ to
each other. The proof is similar to the proof of Case 2.

\begin{figure}[htb]
\begin{center}
\epsfxsize=2.75in \epsfbox{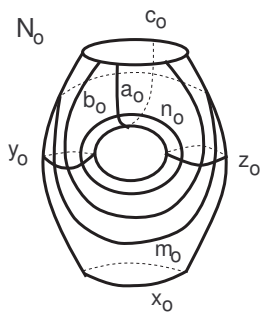} \caption{Arcs and
their encoding circles, II} \label{Figure17}
\end{center}
\end{figure}

Case 5: Assume that $a, b$ are linked, connecting $\partial_i$ to
itself for some $i=1,2.$ By Lemma \ref{D}, there is a
homeomorphism $\phi:(N_o, a_o, b_o) \rightarrow (N, a, b)$ where
$N$ is a regular neighborhood of $\partial_i \cup a \cup b$ in
$S_c$ and $N_o, a_o, b_o$ are as in Figure \ref{Figure17}. Since
$y_o, z_o, c_o$ and $m_o, n_o, c_o$ are the boundary components of
regular neighborhoods of $a_o \cup c_o$ and $b_o \cup c_o$ on
$N_o$ respectively, $\phi(y_o), \phi(z_o)$ and $\phi(m_o),
\phi(n_o)$ are encoding circles for $a$ and $b$ on $S_c$
respectively. We have $(N, \phi(c_o), \phi(x_o), \phi(y_o),
\phi(z_o), \phi(m_o), \phi(n_o))$ $ \cong$ $(N_o, c_o, x_o, y_o,
z_o, m_o, n_o)$. $q(\phi(c_o))$ is equal to $c$ on $S$. So,
$q(\phi(c_o))$ is an essential circle on $S$. Since $c$ is a
nonseparating circle on $S$, $q(\phi(x_o))$ is an essential circle
on $S$. Since $q(\phi(y_o)), q(\phi(z_o))$ and $q(\phi(m_o)),
q(\phi(n_o))$ are encoding circles for essential arcs $a$ and $b$
on $S$ respectively, $q(\phi(y_o)), q(\phi(z_o)), q(\phi(m_o))$
and $q(\phi(n_o))$ are essential circles on $S$. Note that $d \in
\lambda([q(\phi(c_o))])$. Then, by using Lemma \ref{3}, we can
choose $x_o' \in \lambda ([q(\phi(x_o))])$, $y_o' \in
\lambda([q(\phi(y_o))]), z_o' \in \lambda([q(\phi(c_2))]), m_o'
\in \lambda ([q(\phi(m_o))]), n_o' \in \lambda ([q(\phi(n_o))])$,
$N' \subseteq S$ and a homeomorphism $\chi : (N_o, c_o, x_o, y_o,
z_o, m_o$, $n_o) \rightarrow (N', d, x_o', y_o', z_o', m_o',
n_o')$. Since $q(\phi(y_o)), q(\phi(z_o))$ and $q(\phi(m_o)),
q(\phi(n_o))$ are encoding circles for $a$ and $b$ on $S$
respectively, $y_o', z_o'$ and $m_o', n_o'$ are encoding circles
for $\lambda_*([a]), \lambda_*([b])$ on $S$ respectively.
Existence of $\chi$ shows that $\lambda_*([a]), \lambda_*([b])$
have disjoint representatives. $\chi(a_o)$ and $\chi(b_o)$ are
disjoint representatives for $\lambda_*([a]), \lambda_*([b])$
respectively.

We have shown in all the cases that if two vertices have disjoint
representatives, then $\lambda_*$ sends them to two vertices which
have disjoint representatives. Hence, $\lambda_*$ extends to a
simplicial map $\lambda_*: \mathcal{B}(S_c) \rightarrow
\mathcal{B}(S_d)$.\end{proof}

\begin{lemma}
\label{inj2} Let $\lambda : \mathcal{C}(S) \rightarrow
\mathcal{C}(S)$ be a superinjective simplicial map. Then,
$\lambda_* : \mathcal{B}(S_c) \rightarrow \mathcal{B}(S_d)$ is
injective.\end{lemma}

\begin{proof} It is enough to prove that $\lambda_*$ is injective
on the vertex set, $\mathcal{V}(S_c)$. Let $[i], [j] \in
\mathcal{V}(S_c)$ such that $\lambda_*([i])= \lambda_*([j])=[k]$.
Then, by the definition of $\lambda_*$, the type of [i] and [j]
are the same. Assume they are both type 1.1. Let $\{[c], [x],
[y]\}$ and $\{[c], [z], [t]\}$ be the encoding simplices for $[i]$
and $[j]$ respectively. Then, $\{\lambda([c]), \lambda([x]),
\lambda([y])\}$ and $\{ \lambda(c]), \lambda([z]), \lambda([t])\}$
are encoding simplices of $[k]$. So, \{$\lambda ([x]), \lambda
([y])\} = \{\lambda ([z]), \lambda ([t])\}$. Then, since $\lambda$
is injective, we get $\{[x], [y]\} = \{[z], [t]\}$. This implies
$[i]=[j]$. The cases where $[i],[j]$ are both type 1.2 or type 2
can be proven similarly to the first case by using the injectivity
of $\lambda$.\end{proof}

\begin{lemma}
\label{fixing} If an injection $\mu : \mathcal{B}(S_c) \rightarrow
\mathcal{B}(S_d)$ agrees with $h_* : \mathcal{B}(S_c) \rightarrow
\mathcal{B}(S_d)$ on a top dimensional simplex, where $h_*$ is
induced by a homeomorphism $h: S_c \rightarrow S_d$, then $\mu$
agrees with $h_*$ on $\mathcal{B}(S_c)$.
\end{lemma}

\begin{proof} Suppose that $\mu$ agrees with $h_*$ on a top
dimensional simplex, $\Delta$. Then, $h_* ^{-1} \circ \mu$ fixes
$\Delta$ pointwise. Let $w$ be a vertex in $\mathcal{B}(S_c)$. If
$w \in \Delta$, then $w$ is fixed by $h_* ^{-1} \circ \mu$.
Suppose that $w$ is not in $\Delta$. Take a top dimensional
simplex containing $w$. Call it $\Delta'$. There exists a chain
$\Delta = \Delta_0, \Delta_1, ..., \Delta_m = \Delta'$ of top
dimensional simplices in $\mathcal{B}(S_c)$, connecting $\Delta$
to $\Delta'$ such that any two consecutive simplices $\Delta_i,
\Delta_{i+1}$ have exactly one common face of codimension 1. This
follows from the Connectivity Theorem for Elementary Moves of
Mosher, \cite{Mos}, appropriately restated for surfaces with
boundaries. Since $h_* ^{-1} \circ \mu$ is injective, $\Delta_1$
must be sent to a top dimensional simplex by $h_* ^{-1} \circ
\mu$. Let $w_1$ be the vertex of $\Delta_1$ which is not in
$\Delta$. Since $\Delta$ is fixed by $h_* ^{-1} \circ \mu$, the
common face of $\Delta$ and $\Delta_1$ is fixed. $\Delta_1$ is the
unique top dimensional simplex containing the common codimension 1
face of $\Delta$ and $\Delta_1$ other than $\Delta$. Since $h_*
^{-1} \circ \mu (\Delta_1)$ is a top dimensional simplex having
this common face which is different from $\Delta$, $\Delta_1$ must
be sent onto itself by $h_* ^{-1} \circ \mu$. This implies that
$w_1$ is fixed. By an inductive argument, we can prove that all
the top dimensional simplices in the chain are fixed. This shows
that $w$ is fixed. Hence, $h_* ^{-1} \circ \mu$ is the identity on
$\mathcal{B}(S_c)$ and $\mu$ agrees with $h_*$ on
$\mathcal{B}(S_c)$.\end{proof}

We have proven that $\lambda$ is an injective simplicial map which
preserves the geometric intersection 0 and 1 properties. Using
these properties and following N.V.Ivanov's proof of his Theorem
1.1 \cite{Iv1}, it can be seen that $\lambda_*$ agrees with a map,
$h_*$, induced by a homeomorphism $h: S_c \rightarrow S_d$ on a
top dimensional simplex in $\mathcal{B}(S_c)$. Then, by Lemma
\ref{fixing}, it agrees with $h_*$ on $\mathcal{B}(S_c)$. Then,
again by Ivanov's proof, $\lambda$ agrees with a map, $g_*$, which
is induced by a homeomorphism $g: S \rightarrow S$ on
$\mathcal{C}(S)$. This proves Theorem 1.4.

{\bf Remark}: Note that at the end of the proof of Theorem 1.4 we
appealed to Ivanov's proof. In his proof of Theorem 1.1
\cite{Iv1}, by using Lemma 1 in his paper and using some homotopy
theoretic results about the complex of curves on surfaces with
boundary, he shows that an automorphism of the complex of curves
preserves the geometric intersection 1 property if the surface has
genus at least 2. He uses this property to induce automorphisms on
the complex of arcs using automorphisms of the complex of curves
and gets an element of the extended mapping class group which
agrees with the automorphism of the complex of curves. In this
paper, for closed surfaces of genus at least 3, we prove that a
superinjective simplicial map $\lambda$ preserves the geometric
intersection 1 property by using Ivanov's Lemma 1 and using some
elementary surface topology arguments avoiding any homotopy
theoretic results about the complex of curves. These elementary
surface topology arguments can be used to replace the usage of the
homotopy theoretic results about the complex of curves in Ivanov's
proof to show that an automorphism of the complex of curves
preserves the geometric intersection 1 property. Then we follow
Ivanov's ideas and use this property to get an injective
simplicial map of the complex of arcs of $S_c$ (c is
nonseparating) using a superinjective simplicial map of the
complex of curves of $S$ and follow Ivanov to get an element of
the extended mapping class group that we want.

\section{Injective Homomorphisms of Subgroups of Mapping Class Groups}

In this section we assume that $\Gamma'=ker(\varphi)$, where
$\varphi: Mod_S^* \rightarrow Aut(H_1(S, \mathbb{Z}_3))$ is the
homomorphism defined by the action of homeomorphisms on the
homology.

For a group $G$ and for subsets $A, H \subseteq G$, where $A
\subseteq H$, the centralizer of $A$ in $H$, $C_H(A)$, and center
of $G$, $C(G)$, are defined as follows;

 $C_H(A) = \{ h \in H : ha = ah \;\forall a \in A \}$,
 $C(G) = \{ z \in G : zg = gz \;\forall g \in G \}$.

\begin{lemma}
\label{subset}
Let $H$ be a subgroup of a group $G$ and let $A \subseteq H$. Then
  \[
  H \cap C(C_{G} (A)) \subseteq C(C_{H} (A) ).
  \]
\end{lemma}

\begin{proof} Let $h \in  H \cap C(C_{G} (A))$. Then, $h \in H
\cap C_{G}(A)$. It is clear that $H \cap C_{G}(A) = C_{H} (A)$.
So, $h \in C_{H}(A)$. Let $k \in C_{H}(A)$. Then, $k \in C_{G}(A)
\cap H \subseteq C_{G}(A)$. Then, since $h \in C(C_{G} (A))$, $h$
commutes with $k$. Therefore, $h \in C(C_{H} (A))$.\end{proof}

\begin{lemma}
\label{ranklemma} Let $K$ be a finite index subgroup of $Mod_S^*$.
Let $f:K \rightarrow Mod_S^*$ be an injective homomorphism. Let
$\Gamma = f^{-1}(\Gamma') \cap \Gamma'$. Let $G$ be a free abelian
subgroup of $\Gamma$ of rank $3g - 3$. If $a \in G$, then
  \[ \rank C(C_{\Gamma'} (f(a))) \leq \rank C(C_{\Gamma}(a)). \]
\end{lemma}

\begin{proof} Let $A= f(G) \cap C(C_{\Gamma'} (f(a)))$ and $B =
\langle f(G), C(C_{\Gamma'} (f(a))) \rangle$ be the group
generated by $f(G)$ and $C(C_{\Gamma'} (f(a)) )$. Since $a \in G$
and $G$ is abelian, $f(a) \in f(G)$ and $f(G)$ is abelian. Then,
$f(G) \subseteq C_{\Gamma'}(f(a))$ since $f(G) \subseteq \Gamma'$.
Then, B is an abelian group. Since rank $ G= 3g-3$ and $f$ is
injective on G, rank $f(G) = 3g-3$. The maximal rank of an abelian
subgroup of $Mod_S^*$ is $3g-3$, \cite{BLM}. So, since $f(G)
\subseteq B$ and $B$ is an abelian group and rank $f(G) = 3g-3$ we
have rank $B = 3g-3$.

 We have the following exact sequence,
 \[ 1 \rightarrow A \rightarrow f(G) \oplus  C(C_{\Gamma'} (f(a)) )
 \rightarrow B \rightarrow 1 \]

  This gives us
  \[ \rank f(G) + \rank C(C_{\Gamma'} (f(a))) = \rank A + \rank B \]

 Since $\rank f(G) = \rank B$, we get
 \begin{equation}
   \label{rankequation}
   \rank C(C_{\Gamma'} (f(a))) = \rank A.
 \end{equation}

Then, by using Lemma \ref{subset}, we get $A= f(G) \cap
C(C_{\Gamma'}(f(a))) \subseteq$ $f(\Gamma) \cap C(C_{\Gamma'}
(f(a)))$ $\subseteq C(C_{f(\Gamma)} (f(a))). $ Since $f$ is
injective on $K$, $C(C_{f(\Gamma)} (f(a)))$ is isomorphic to
$C(C_{\Gamma}(a))$. So, we get rank $A \leq$ rank $C(C_{\Gamma}
(a))$. Then, by equation \ref{rankequation}, we have
 \[
        \rank C(C_{\Gamma'} (f(a))) \leq  \rank C(C_{\Gamma}(a))
 \]  \end{proof}

%A homeomorphism $f: S \rightarrow S$ is called {\it pure} if for
%some closed one-dimensional submanifold (possibly empty) $C$ of
%$S$ the following condition is satisfied: The components of $C$
%are nontrivial (that is, they do not deform on $S$ to a point) and
%are pairwise nonisotopic, $F$ is fixed (i.e., is the identity) on
%$C$, it does not rearrange the components of $S \setminus C$, and
%it induces on each component of $S_C$ a homeomorphism isotopic to
%either a pseudo-Anosov or the identity homeomorphism. A mapping
%class is called {\it pure} if it has a pure representative.

\begin{lemma}
\label{center1}  Let $\Gamma$ be a finite index subgroup of
$Mod_S^*$ and $\Gamma \subseteq \Gamma'$. Then, $\exists N \in
\mathbb{Z^*}$ s.t. $\forall \alpha \in \mathcal{A}$, $C(C_{\Gamma}
(t_{\alpha} ^{N}))$ is an infinite cyclic subgroup of $\langle
t_\alpha \rangle$.
\end{lemma}

\begin{proof} Let $N \in \mathbb{Z}^*$ such that $t_\alpha ^N \in
\Gamma$ for all $\alpha \in \mathcal{A}$. Let $\alpha, \beta \in
\mathcal{A}$ such that $i(\alpha, \beta)= 0$. Let $f \in
C(C_{\Gamma} (t_{\alpha} ^{N}) )$. Then, $f(\alpha)= \alpha$. We
want to show that $f(\beta)= \beta$. Since $i(\alpha, \beta)=0$,
$t_\alpha$ commutes with $t_\beta$. Then, $t_{\alpha} ^N$ commutes
with $t_\beta^N$. Then, since $t_\beta^N \in \Gamma$, $t_\beta^N
\in C_{\Gamma} (t_{\alpha} ^{N})$. Since $f \in C(C_{\Gamma}
(t_{\alpha} ^{N}))$, $f$ commutes with $t_\beta^N$. Then,
$f(\beta)= \beta$. Since $f$ fixes the isotopy class of every
circle which has geometric intersection zero with $\alpha$, the
reduction, $\hat{f}$, of $f$ along $\alpha$ (\cite{BLM}) fixes the
isotopy class of every circle on $S_a$ where $a \in [\alpha]$.
Since $\Gamma \subseteq \Gamma'$, $f \in \Gamma'$. Hence, by Lemma
1.6 \cite{Iv2}, $\hat{f}$ restricts to each component $Q$ of $S_a$
and the restriction of $\hat{f}$ to $Q$ is either trivial or
infinite order. By Lemma 5.1 \cite{IMc} and Lemma 5.2 \cite{IMc},
this restriction is finite order. Hence, it is trivial. Then,
$f=t_\alpha ^r$ for some $r \in \mathbb{Z}$, \cite{BLM}. So,
$C(C_{\Gamma} (t_{\alpha} ^{N}) ) \subseteq$ $\langle t_\alpha
\rangle$. Since $t_\alpha ^{N} \in C(C_{\Gamma} (t_{\alpha}
^{N}))$, $C(C_{\Gamma} (t_{\alpha} ^{N}) )$ is a nontrivial
subgroup of an infinite cyclic group. Hence, it is infinite
cyclic. \end{proof}

\begin{lemma}
\label{rank=1} Let $K$ be a finite index subgroup of $Mod_S^*$ and
$f:K \rightarrow Mod_S^*$ be an injective homomorphism. Let
$\alpha \in \mathcal{A}$. Then there exists $N \in \mathbb{Z^*}$
such that

\begin{center}
 $rank$ $C(C_{\Gamma'} (f(t_{\alpha} ^{N})) ) = 1$.
\end{center}
\end{lemma}

\begin{proof} Let $\Gamma=f^{-1}(\Gamma') \cap \Gamma'$. Let $\alpha \in
\mathcal{A}$, $N \in \mathbb{Z^*}$ such that $t_{\alpha} ^ N \in
\Gamma$. Since $\langle f(t_\alpha ^N) \rangle \subseteq
C(C_{\Gamma'}(f(t_{\alpha}^{N})))$, and $\langle f(t_\alpha
^N)\rangle$ is an infinite cyclic group we have
\begin{equation}
\label{1}
 \rank C(C_{\Gamma'}(f(t_{\alpha} ^{N}))) \geq 1.
\end{equation}

By Lemma \ref{center1}, $\rank C(C_{\Gamma} (t_{\alpha}^{N})) = 1$. Since $t_{\alpha}^{N}$ is in a
free abelian subgroup of $\Gamma$ of rank $3g-3$, by Lemma \ref{ranklemma}, $\rank C(C_{\Gamma'}
(f(t_\alpha ^{N}))) \leq \rank C(C_{\Gamma}(t_\alpha^{N}))$. So, we get
\begin{equation}
\label{2}
 \rank C(C_{\Gamma'} (f(t_\alpha ^{N}))) \leq 1.
\end{equation}

(\ref{1}) and (\ref{2}) imply that $ \rank C(C_{\Gamma'}(f(t_{\alpha} ^N))) = 1$. \end{proof}

\begin{lemma}
\label{reducible} Let K be a finite index subgroup of $Mod_S^*$.
Let $f:K \rightarrow Mod_S^*$ be an injective homomorphism. Then
there exists $N \in \mathbb{Z^*}$ such that $f(t_{\alpha} ^ N)$ is
a reducible element of infinite order for all $\alpha \in
\mathcal{A}$.
\end{lemma}

\begin{proof} Let $N \in \mathbb{Z^*}$ such that
$t_{\alpha} ^ N \in K$ for all $\alpha \in \mathcal{A}$. Since $f$
is an injective homomorphism on $K$ and $t_{\alpha} ^ N$ is an
infinite order element of $K$, $f(t_{\alpha} ^ N)$ is an infinite
order element of $Mod_S^*$. So, $f(t_{\alpha} ^ N)$ is either a
reducible element or p-Anosov. Suppose it is p-Anosov. Let $J$ be
a maximal system of circles containing $a$ where $[a]=\alpha$ and
$J'$ be the set of isotopy classes of these circles. Let $T_{J'}$
be the subgroup of $Mod_S^*$ generated by $t_{\beta} ^ N$ for all
$\beta$ in $J'$. $T_{J'}$ is a free abelian subgroup of K and it
has rank $3g-3$. Since $f$ is an injection, $f(T_{J'})$ is also a
free abelian subgroup of rank $3g-3$. It contains $f(t_{\alpha} ^
N)$. By \cite{Mc}, $C_{Mod_S^*}(f(t_{\alpha} ^ N))$ is a virtually
infinite cyclic group and $f(T_{J'}) \subseteq
C_{Mod_S^*}(f(t_{\alpha} ^ N))$. Then $3g-3 \leq 1$. This gives us
$3g \leq 4$. But this is a contradiction to the assumption that $g
\geq 3$. Hence, $f(t_{\alpha} ^ N)$ is a reducible element of
infinite order.
\end{proof}

\begin{lemma}
\label{wellknown1} Let $\alpha, \beta \in \mathcal{A}$ and $i, j$ be nonzero integers. Then, $t_\alpha ^i
=t_\beta ^j \Leftrightarrow \alpha = \beta$ and $i=j$.
\end{lemma}

\begin{proof} If $\alpha = \beta$ and $i=j$, then obviously $t_\alpha ^i = t_\beta ^j$. To see the other
implication, we consider the reduction systems. Since $\beta$ is the canonical reduction system for $t_\beta^j$
and $\alpha$ is the canonical reduction system for $t_\alpha^i$ and $t_\alpha ^i= t_\beta ^j$, their canonical
reduction systems must be equal. So, $\alpha=\beta$. Then we have $t_\alpha ^i = t_\alpha ^j$. Then, $i=j$
since $t_\alpha$ is an infinite order element in $Mod_S^*$.  \end{proof}

\begin{lemma}
\label{correspondence} Let $K$ be a finite index subgroup of
$Mod_S^*$, and $f:K \rightarrow Mod_S^*$ be an injective
homomorphism. Then $\forall \alpha \in \mathcal{A}$, $f( t_\alpha
^N)= t_{\beta(\alpha)}^M$ for some $M, N \in \mathbb{Z^*}$,
$\beta(\alpha) \in \mathcal{A}$.
\end{lemma}

\begin{proof} Let $\Gamma= f^{-1}(\Gamma') \cap \Gamma'$. Since
$\Gamma$ is a finite index subgroup we can choose $N \in Z^*$ such
that $t_\alpha^N \in \Gamma$, for all $\alpha$ in $\mathcal{A}$.
By Lemma \ref{reducible}, $f(t_{\alpha} ^ N)$ is a reducible
element of infinite order in $Mod_S^*$. Let $C$ be a realization
of the canonical reduction system of $f(t_{\alpha}^N)$. Let $c$ be
the number of components of $C$ and $p$ be the number of p-Anosov
components of $f(t_{\alpha} ^N)$. Since $t_{\alpha} ^ N \in
\Gamma, f(t_{\alpha} ^ N) \in \Gamma'$. By Theorem 5.9 \cite{IMc},
$C(C_{\Gamma'} (f(t_{\alpha} ^ N )))$ is a free abelian group of
rank $c+p$. By Lemma \ref{rank=1}, $c+p=1$. Then, either $c=1$,
$p=0$ or $c=0$, $p=1$. Since there is at least one curve in the
canonical reduction system we have $c=1$, $p=0$. Hence, since
$f(t_{\alpha} ^ N) \in \Gamma'$, $f(t_{\alpha} ^{N}) = t_{\beta
({\alpha})}^{M}$ for some $M \in \mathbb{Z^*}$, $\beta(\alpha) \in
\mathcal{A}$, [1], [3].\end{proof}

{\bf Remark}: Suppose that $f(t_{\alpha} ^{M}) = t_\beta ^P$ for
some $\beta \in \mathcal{A}$ and $M, P \in \mathbb{Z^*}$ and
$f(t_{\alpha} ^{N}) = t_\gamma ^Q$ for some $\gamma \in
\mathcal{A}$ and $N, Q \in \mathbb{Z^*}$. Since $f(t_{\alpha} ^{M
\cdot N}) = f(t_{\alpha} ^{N \cdot M})$, $t_\beta ^{PN} = t_\gamma
^{QM}$, $P, Q, M, N \in \mathbb{Z^*}$. Then, $\beta = \gamma$ by
Lemma \ref{wellknown1}. Therefore, by Lemma \ref{correspondence},
$f$ gives a correspondence between isotopy classes of circles and
$f$ induces a map, $f_*: \mathcal{A} \rightarrow \mathcal{A}$,
where $f_*(\alpha) = \beta(\alpha)$.

In the following lemma, we use a well known fact that $f t_\alpha
f^{-1}=t_{f(\alpha)} ^{\epsilon(f)}$ for all $\alpha$ in
$\mathcal{A}$, $f \in Mod_S^*$, where $\epsilon(f) = 1$ if $f$ has
an orientation preserving representative and $\epsilon(f) = -1$ if
$f$ has an orientation reversing representative.

\begin{lemma}
\label{identity} Let $K$ be a finite index subgroup of $Mod_S^*$.
Let $f:K \rightarrow Mod_S^*$ be an injective homomorphism. Assume
that there exists $N \in \mathbb{Z}^*$ such that $\forall \alpha
\in$ $\mathcal{A}$, $\exists Q \in \mathbb{Z}^*$ such that
$f(t_{\alpha} ^N) = t_{\alpha}^Q$. Then, $f$ is the identity on
$K$.
\end{lemma}

\begin{proof}
We use Ivanov's trick to see that $f(kt_{\alpha} ^ N k^{-1})=$
$f(t_{k(\alpha)} ^{\epsilon{(k)} \cdot N }) = t_{k(\alpha)} ^{Q
\cdot \epsilon{(k)}}$ and $f(kt_{\alpha} ^ N k^{-1}) = f(k)
f(t_{\alpha} ^N) f(k)^{-1}=$ $f(k) t_{\alpha} ^Q f(k)^{-1} =
t_{f(k)(\alpha)} ^{\epsilon(f(k))\cdot Q}$ $\forall \alpha \in
\mathcal{A}$, $\forall k \in K$. Then, we have $t_{k(\alpha)} ^{Q
\cdot \epsilon{(k)}} = t_{f(k)(\alpha)} ^{\epsilon(f(k)) \cdot Q}
\Rightarrow$ $k(\alpha) = f(k)(\alpha)$ $\forall \alpha \in
\mathcal{A}$, $\forall k \in K$ by Lemma \ref{wellknown1}. Then,
$k^{-1}f(k)(\alpha) = \alpha$ $\forall \alpha \in \mathcal{A}$,
$\forall k \in K$. Then, $k^{-1}f(k)$ commutes with $t_{\alpha}$
$\forall \alpha \in \mathcal{A}$, $\forall k \in K$. Since
$Mod_S^*$ is generated by Dehn twists, $k^{-1}f(k) \in C(Mod_S^*)$
$\forall k \in K$. Since the genus of $S$ is at least 3,
$C(Mod_S^*)$ is trivial. So, $k = f(k)$ $\forall k \in K$. Hence,
$f=id_K$.
\end{proof}

\begin{coroll}
\label{id} Let $g: Mod_S^* \rightarrow Mod_S^*$ be an isomorphism
and $h : Mod_S^* \rightarrow Mod_S^*$ be an injective
homomorphism.  Assume that there exists $N \in \mathbb{Z}^*$ such
that $\forall \alpha \in$ $\mathcal{A}$, $\exists Q \in
\mathbb{Z}^*$ such that $h(t_{\alpha} ^N) = g(t_{\alpha}^Q)$. Then
$g=h$.
\end{coroll}

\begin{proof} Apply Lemma \ref{identity} to $g^{-1} h$ with $K = Mod_S^*$.
Since for all $\alpha$ in $\mathcal{A}$, $g^{-1} h(t_{\alpha} ^N)
= t_{\alpha} ^{Q}$, we have $g^{-1} h = id_K$. Hence, $g =
h$.\end{proof}

\begin{lemma}
\label{wellknown2} Let $\alpha, \beta$ be distinct elements in $\mathcal{A}$. Let $i, j$ be two nonzero
integers. Then, $t_\alpha ^i t_\beta ^j = t_\beta ^j t_\alpha ^i \Leftrightarrow i(\alpha, \beta)= 0$.
\end{lemma}

\begin{proof} $t_\alpha ^i t_\beta ^j = t_\beta ^j t_\alpha ^i \Leftrightarrow
t_ \alpha ^i t_\beta ^j t_\alpha ^{-i} = t_\beta ^j
\Leftrightarrow t_{t_\alpha^i (\beta) } ^j  = t_\beta ^j
\Leftrightarrow t_\alpha ^i(\beta) = \beta$ (by Lemma
\ref{wellknown1}) $\Leftrightarrow i(\alpha, \beta)=0$
(``$\Leftarrow$" : Clear. ``$\Rightarrow$": $i(\alpha, \beta) \neq
0$ $\Rightarrow$ $t_\alpha ^m(\beta) \neq \beta$ for all $m \in
\mathcal{Z}^*$ since $\alpha$ is an essential reduction class for
$t_\alpha$.). \end{proof}

By the remark after Lemma \ref{correspondence}, we have that $f: K
\rightarrow Mod_S^*$ induces a map $f_*: \mathcal{A} \rightarrow
\mathcal{A}$. We prove the following lemma which shows that $f_*$
is a superinjective simplicial map on $\mathcal{C}(S)$.

\begin{lemma}
\label{intersection0} Let $f:K \rightarrow Mod_S^*$ be an injection.
Let $\alpha$, $\beta \in \mathcal{A}$. Then,
\[
i(\alpha,\beta)=0 \Leftrightarrow i(f_{*}(\alpha),
f_{*}(\beta))=0.
\]
\end{lemma}

\begin{proof} There exists $N \in \mathbb {Z^*}$ such that
$t_{\alpha} ^N \in K$ and $t_{\beta} ^N \in K$. Then we have the
following: $i(\alpha, \beta)=0$ $\Leftrightarrow$ $t_{\alpha} ^N
t_\beta ^N = t_{\beta} ^N t_\alpha ^N$ (by Lemma \ref{wellknown2})
$\Leftrightarrow$ $f(t_{\alpha} ^N) f(t_{\beta} ^ N) = f(t_{\beta}
^ N) f(t_{\alpha} ^ N)$ (since $f$ is injective on K)
$\Leftrightarrow$ $t_{f_*(\alpha)} ^P t_{f_*(\beta)}^Q =
t_{f_*(\beta)}^Q t_{f_*(\alpha)} ^P$ where $P = M(\alpha_{1}, N),
Q = M(\alpha_{2}, N) \in \mathbb{Z}^*$ (by Lemma
\ref{correspondence}) $ \Leftrightarrow i(f_{*}(\alpha),
f_{*}(\beta))=0$ (by Lemma \ref{wellknown2}).\end{proof}

Now, we prove the second main theorem of the paper.

\begin{theorem}
\label{main} Let $f$ be an injective homomorphism, $f:K
\rightarrow Mod_S^*$, then $f$ is induced by a homeomorphism of
the surface $S$ and $f$ has a unique extension to an automorphism
of $Mod_S^*$. \end{theorem}

\begin{proof} By Lemma \ref{intersection0}, $f_*$ is a superinjective
simplicial map on $\mathcal{C}(S)$. Then, by Theorem 1.4, $f_*$ is
induced by a homeomorphism $h:S \rightarrow S$, i.e. $f_*(\alpha)
= H(\alpha)$ for all $\alpha$ in $\mathcal{A}$, where $H=[h]$. Let
$\chi ^ {H}: Mod_S^* \rightarrow Mod_S^*$ be the isomorphism
defined by the rule $\chi ^ {H}(K) = HKH^{-1}$ for all $K$ in
$Mod_S^*$. Then for all $\alpha$ in $\mathcal{A}$, we have the
following:

$\chi ^{H ^{-1}} \circ f ({t_ \alpha} ^N) =  \chi ^{H ^{-1}}
(t_{f_*(\alpha)}^ M) = \chi ^{H ^{-1}} (t_{H(\alpha)} ^M) = H^{-1}
t_{H(\alpha)} ^M H = t_ {H^{-1} (H(\alpha))} ^{M \cdot
\epsilon{(H^{-1})}} = t_\alpha ^{M \cdot \epsilon{(H^{-1})}}$.

Then, since $\chi ^{H^{-1}} \circ f$ is injective, $\chi ^{H^{-1}}
\circ f = id_K$ by Lemma \ref{identity}. So, $\chi ^H |_K = f$.
Hence, $f$ is the restriction of an isomorphism which is
conjugation by $H$, (i.e. $f$ is induced by $h$).

Suppose that there exists an automorphism $\tau :  Mod_S^*
\rightarrow Mod_S^*$ such that $\tau |_{K}=f$. Let $N \in Z^*$
such that $ t_\alpha ^N \in K$  for all $\alpha$ in $\mathcal{A}$.
Since $\chi ^H |_K = f = \tau |_K$ and $t_\alpha ^N \in K$,
$\tau(t_\alpha ^N) = \chi ^H(t_\alpha ^N)$ for all $\alpha$ in
$\mathcal{A}$. Then, by Corollary \ref{id}, $\tau = \chi ^ {H}$.
Hence, the extension of $f$ is unique. \end{proof}\\

{\bf Acknowledgments}

This work is the author's doctoral dissertation at Michigan State
University. We thank John D. McCarthy, the author's thesis
advisor, for his guidance throughout this work and for reviewing
this paper. We also thank Nikolai Ivanov for his suggestions about
this work.

%\newpage
%\listoffigures

\end{document}